\documentclass[final,hidelinks,onefignum,onetabnum]{siamart251216}

\usepackage[utf8]{inputenc}
\usepackage{amssymb,amsfonts}
\usepackage{graphicx}

\newcommand{\E}{\mathbb{E}}
\newcommand{\Epath}{\mathbb{E}_{\mathrm{path}}}
\newcommand{\Erf}{\mathbb{E}_{\mathrm{rf}}}
\newcommand{\R}{\mathbb{R}}
\newcommand{\Ac}{\mathcal{A}}
\newcommand{\Xc}{\mathcal{X}_{\mathrm{in}}}
\newcommand{\Ec}{\mathcal{E}}

\newsiamthm{assumption}{Assumption}
\newtheorem*{remark}{Remark}

\title{Random Feature Neural Networks for Stochastic Optimal Control}
\author{Guangyi He\thanks{Department of Mathematics, Imperial College London, London, UK.} \and Damir Filipovic\thanks{Swiss Finance Institute, \'Ecole Polytechnique F\'ed\'erale de Lausanne, Lausanne, Switzerland.} \and Lukas Gonon\thanks{School of Computer Science, University of St.\ Gallen, St.\ Gallen, Switzerland, and Department of Mathematics, Imperial College London, London, UK.}}
\date{\today}

\begin{document}

\maketitle

\begin{abstract}
Random feature neural networks (RFNNs) sample their hidden parameters once and train only the output layer. We analyze RFNNs as feedback policies for finite-horizon stochastic control with trajectory-dependent objectives. A policy change alters later states, so the error that matters is the objective error of the induced trajectory rather than a pointwise regression error. We extend the quantitative random-feature approximation theory for PDE solutions and nonlinear regression to feedback control. We investigate the objective error of an RFNN policy through two sources, approximation of the policy by random features and estimation of the objective from finitely many paths. A stability estimate for the controlled trajectory turns a policy error into an objective error. Under bounded actions, a regular benchmark policy and exact empirical minimization, the expected objective gap to the benchmark policy is at most of order $N^{-1/2}+\sqrt{N(1+\log N)/n}$, with $N$ the feature width per output coordinate and $n$ the number of paths. At fixed sample size the two terms trade off in $N$, and balancing them gives $N\asymp\sqrt{n/\log n}$.  Across a wide range of examples, estimated gaps generally decrease with width and sample size, and RFNNs show an advantage when simulated paths are limited. We also study Parsimonious Random Feature Learning (PRFL), which grows and prunes the random dictionary during training while training only the output coefficients. At matched retained width and training budget, PRFL improves on a fixed RFNN in the LQ and Heston comparisons, while the fully trained network remains best in those examples. 
\end{abstract}

\section{Introduction}

In finite-horizon stochastic control problems, optimal policies are often unavailable in closed form. A common approach is to parameterize the policy and optimize its objective using simulated trajectories. The choice of policy representation is central to this approach, because a richer representation can express better policies but is also more prone to overfitting when training data are limited. Neural networks provide flexible policy representations for high-dimensional and path-dependent problems \cite{bachouch2022deep,buehler2019deep}, but understanding their approximation and statistical properties remains important for assessing the quality of learned policies.

In this paper, we study random feature neural networks (RFNNs), which fix the hidden parameters at random values and train only the output layer. This structure reduces the number of trainable parameters and makes the policy linear in the output weights. The current literature on random feature models mainly provides quantitative approximation and learning results on regression and PDE problems \cite{Gonon2023random,lanthaler2023vector}. In feedback control, a policy change alters subsequent states and decisions, so its effect on the objective accumulates over the trajectory.

We study these questions in a common framework for Markovian control and history-dependent policy optimization. The state $X_t$ evolves under a transition that combines the current state, the action and a fresh innovation, and a feedback policy maps the state and the decision date to the action. The state may contain the observed history or an appropriate problem-specific summary, including quantities determined by previous actions. The problem is to choose a feedback policy $H$ to minimize
\begin{equation}
\label{eq:draft_objective}
J(H)=\inf_{\beta\in\mathcal B}\Epath\left[\mathcal L(\beta,X^H,\alpha^H)\right],
\end{equation}
where $X^H$ is the full state trajectory induced by $H$, $\alpha^H_t=H(X_t^H,t)$, $\mathcal L$ is the trajectory loss, and $\beta$ is an auxiliary optimization variable. This form includes standard expected costs and risk objectives with an optimized certainty equivalent representation. Together with the flexible state representation, this provides a common formulation for both Markovian control and path-dependent problems such as hedging with transaction costs.

We investigate the objective error by splitting it into an approximation error and a statistical error. The first is the approximation error of the RFNN class itself. Gonon's scalar uniform bound controls the policy approximation error at order $O(N^{-1/2})$ in expectation over the features \cite[Theorem~1]{Gonon2023random}, where $N$ is the feature width per output coordinate. A stability estimate for the controlled trajectory propagates this to the overall performance, so the objective retains the $N^{-1/2}$ rate, with a constant that grows polynomially in the dimension whenever the benchmark-policy constants do. The second is the statistical error, the discrepancy between the empirical and population objectives computed from $n$ independent simulated paths. A uniform empirical-process bound over the policy class gives order $\sqrt{N(1+\log N)/n}$, which vanishes as $n$ grows but increases with $N$. Combining the two yields
\begin{equation}
\label{eq:draft_excess_bound}
\Erf\,\Epath\bigl[J(\widehat H)-J(H^\star)\bigr]
\leq O\!\left(\frac{1}{\sqrt N}+\sqrt{\frac{N(1+\log N)}{n}}\right).
\end{equation}
At fixed sample size, a larger width lowers the approximation error but raises the statistical one, and balancing them gives $N\asymp\sqrt{n/\log n}$ and an upper bound of order $((\log n)/n)^{1/4}$. 

We evaluate the analysis on a wide range of examples. The Markovian control examples, linear-quadratic and semilinear, follow \cite{bachouch2022deep}, and the Heston deep hedging example follows \cite{buehler2019deep}. We estimate the approximation error and the statistical error, which decrease in $N$ and $n$ respectively as expected. At fixed sample size, a wider representation lowers the approximation error but raises the statistical contribution, showing the trade-off between the two errors. In the comparison between RFNN and fully trained neural networks (FTNNs), RFNNs show an advantage when simulated paths are limited.

We also investigate adaptive feature selection through Parsimonious Random Feature Learning (PRFL). Drawing on parsimonious kernel and random-feature methods \cite{koppel2019parsimonious,lemoal2025parsimonious}, PRFL starts from a small random dictionary and repeatedly trains the output layer, removes features with small squared output weights and adds a fresh block of features. This keeps the output-layer optimization problem of the RFNN analysis while allowing the policy to draw on a larger feature pool. The trajectory-stability estimate bounds the objective change caused by each pruning step in terms of the squared-weight pruning budget. At matched retained width and training budget, PRFL improves on the fixed RFNN in the LQ and Heston comparisons, while the fully trained network remains best in the LQ and Heston examples.

\section{Related Work}
\label{sec:related_work}

\paragraph{Neural networks  and backward dynamic programming}

Neural-network methods for stochastic control  often retain the backward structure of dynamic programming. An overview of neural-network algorithms for stochastic control and PDEs in finance is given in \cite{germain2023survey}. In \cite{bachouch2022deep}, performance-iteration methods learn feedback controls backward in time, while hybrid methods combine control learning with neural regression of the value function. Several classes of examples in that study, including semilinear PDEs and linear-quadratic control, are also covered by our analysis. The companion analysis in \cite{hure2021deep} bounds the value-function errors of these backward schemes in terms of neural-network approximation and statistical estimation errors. It traces how both errors accumulate over the remaining decision dates.

The work of \cite{dellavecchia2025error} also studies approximation inside dynamic programming but replaces neural regression with kernel ridge regression. Its empirical Bellman operators introduce regression and Monte Carlo errors in addition to the error inherited from later decision dates. The resulting analysis describes how these errors propagate backward through time under stability conditions on the Bellman operator, with American option pricing used as the principal application. These  works study the accumulation of local errors in backward dynamic programming. By contrast, our stability analysis traces the effect of policy approximation forward through the controlled trajectory and transfers the resulting state and action errors to the final objective.

\paragraph{Direct policy optimization and deep hedging}

Direct policy methods parameterize the control itself and optimize it using complete simulated or observed trajectories. Deep Learning Approximation for Stochastic Control Problems \cite{han2016deepcontrol} connects time-dependent control networks through the system dynamics and jointly trains them against the total trajectory cost using simulated noise paths. Deep Hedging \cite{buehler2019deep} applies this approach to financial trading under market frictions, constraints, and convex risk measures. The policy may depend on current market information and previous positions, allowing the strategy to account for the cost of rebalancing. 

The direct policy approach has also been developed along several other directions. For mean-field control and games, \cite{carmona2022convergence} analyzes neural-network methods for McKean--Vlasov dynamics and proves an error rate for a computable surrogate of the original control problem. In deep hedging with a single underlying asset and proportional transaction costs, the optimal strategy has a no-trade band under the assumptions of \cite{imaki2021noarbitrage}. That paper builds the band into the network architecture to speed up training. Second-order optimization has also been used to accelerate the training of hedging policies \cite{mueller2024secondorder}. Robust deep hedging studies strategies that remain effective under model uncertainty \cite{lutkebohmert2021robust}, and distributional adversarial training studies how the training procedure responds to perturbations of the market distribution \cite{he2025distributional}.

Backward dynamic programming has also been applied to deep hedging when the risk criterion admits a recursive representation. The actor--critic formulation of \cite{murray2022continuous} relies on a Bellman recursion for the exponential-utility certainty equivalent. At each decision date, this recursion evaluates the current reward together with a learned continuation value under the same risk criterion. The actor and critic updates rely on this recursion, so changing the risk criterion may require a new dynamic formulation. Our objective \eqref{eq:draft_objective} is instead specified through the loss of the complete controlled trajectory, with auxiliary optimization variables included so that the risk objectives above are covered.

\paragraph{Randomized neural-network approximation}

Randomized neural networks have been well-studied in nonlinear regression and PDE approximation. Extreme learning machines use randomly selected hidden nodes and analytically determined output weights \cite{huang2006extreme}. Earlier, \cite{barron1993universal} established approximation bounds for sigmoidal networks under a Fourier first-moment condition and compared them with lower bounds for fixed-basis expansions.  Random features were introduced for the approximation of kernel machines at large scale \cite{rahimi2007random}, and random bases were shown to approximate functions uniformly \cite{rahimi2008uniform}. Quantitative bounds of this type cover functions, functionals and reservoir systems \cite{gonon2023approximation}, and models with Banach-space-valued targets have their own universal approximation results \cite{neufeld2023banach}. For nondegenerate Black--Scholes type PDEs, \cite{Gonon2023random} obtains uniform bounds of order $N^{-1/2}$ in the width, with constants growing at most polynomially in dimension, under specified Fourier, sampling, and regularity assumptions. Derivative estimates extend these bounds to higher-order Sobolev norms for physics-informed extreme learning machines \cite{deryck2025approximation} and to time-dependent mixed-Sobolev norms for nonlinear PDE examples \cite{mehmood2026rann}. Random neural networks have also been used in numerical schemes for nonlinear PDEs and PIDEs. The random deep splitting method of \cite{neufeld2025splitting} comes with a convergence analysis, while \cite{jacquier2026rough} uses a reservoir-type network and least-squares regression for path-dependent equations arising in rough volatility. For optimal stopping, \cite{herrera2024optimal} uses fixed random features and linear readouts to approximate continuation values on exogenous market paths and proves convergence of the resulting methods. In our setting the networks represent feedback controls that affect subsequent states, so we apply Gonon's scalar uniform bound \cite[Theorem~1]{Gonon2023random} coordinatewise, and our stability analysis transfers the $N^{-1/2}$ rate to the control objective, with a constant that grows polynomially in the dimension and in the benchmark-policy constants.

\paragraph{Statistical learning and error decomposition}
Our analysis of policies trained on simulated paths follows the approximation--estimation framework of statistical learning \cite{bottou2008tradeoffs}. For an exact empirical minimizer, this framework bounds the excess population objective by a best-in-class approximation term and a uniform discrepancy between empirical and population objectives. Rademacher-complexity bounds control this discrepancy over the candidate class \cite{BartlettMendelson2002}. 

In random-feature learning, this decomposition links representation size to sample size.  For ridge regression, order $\sqrt n\log n$ random features suffice for an $O(n^{-1/2})$ excess prediction-risk bound under the assumptions of \cite{rudi2017generalization}. Randomized constructions of the output weights have also been studied as an alternative to optimizing them \cite{rahimi2008weighted}. Related approximation and generalization bounds cover vector-valued outputs \cite{lanthaler2023vector}. In these prediction problems, the input distribution is prescribed and does not change with the fitted coefficients.

For feedback policies, the states at which the policy is evaluated can depend on earlier decisions. Using dynamic programming and a mean-field formulation of single-hidden-layer networks, \cite{baros2025meanfield} establishes generalization bounds for regularized controls learned from finitely many observations. A related pathwise approach is used for optimal stopping under model uncertainty in \cite{belomestny2022empirical}, where finite-sample bounds for a variance-penalized martingale-dual criterion depend on the covering complexity of the martingale class.  

\paragraph{Parsimonious online learning with kernels and random features}

Parsimonious Online Learning with Random Features (POLRF) \cite{lemoal2025parsimonious} shares the RFNN architecture of sampled hidden features and a trainable linear readout, but grows and sparsifies the feature set during learning. It adapts the stochastic-gradient and sparse-projection approach of Parsimonious Online Learning with Kernels (POLK) \cite{koppel2019parsimonious} to online regularized regression.  In its stochastic-control applications, these regressions approximate value or continuation functions within backward dynamic programming. It also provides convergence results for online regularized regression under specified learning-rate and sparsification-budget schedules.

Inspired by POLRF, our Parsimonious Random Feature Learning (PRFL) method adapts feature growth and sparsification to direct policy optimization. It trains an enlarged random-feature policy under the full trajectory loss and prunes features according to their cumulative squared output weights.

\section[Setup]{Basic Setup}
\label{sec:setup}

\subsection{Finite Horizon Stochastic Control Problem}
\label{sec:control_problem}

In this section we consider the general form of the problem investigated in this paper. Let $1,\ldots,T+1$ be the dates at which the state is observed, with decisions made at dates $1,\ldots,T$. The state is $X_t=(Y_t,Z_t)\in\mathcal Y\times\mathcal Z$, where $Y_t$ collects the exogenous coordinates, whose evolution does not depend on the policy, and $Z_t$ the coordinates that the actions drive.

Let the action set be the box $\Ac=\prod_{i=1}^r[a_i,b_i]\subseteq\R^r$. The policy input space is $\Xc=\mathcal X\times\{1,\ldots,T\}\subseteq\R^m$, where $m=d_X+1$ and $d_X$ is the dimension of the full state $X_t=(Y_t,Z_t)$. The policy is evaluated at the decision dates $1,\ldots,T$. For measurable maps $\Phi_t:\mathcal Y\times\Ec\to\mathcal Y$ and $\Psi_t:\mathcal Y\times\mathcal Z\times\Ac\times\Ec\to\mathcal Z$, and a measurable feedback policy $H:\Xc\to\Ac$, the state and the action satisfy
\begin{equation}
Y_{t+1}=\Phi_t(Y_t,\eta_t),
\qquad
Z_{t+1}=\Psi_t(Y_t,Z_t,\alpha_t,\eta_t),
\qquad
\alpha_t=H(X_t,t),
\label{eq:controlled_state_recursion}
\end{equation}
for $t=1,\ldots,T$, with the initial state fixed at $(Y_1,Z_1)=X_1$. The policy enters only through the action, and we write $(X^H,\alpha^H)$ for the trajectory induced by $H$. The innovation sequence $\omega=(\eta_1,\ldots,\eta_T)$ has i.i.d.\ entries in $\Ec$.

The trajectory loss $\mathcal L(\beta,X,\alpha)$ assigns a real-valued loss to the complete state trajectory $X=(X_t)_{t=1}^{T+1}$ and action trajectory $\alpha=(\alpha_t)_{t=1}^{T}$, with an optional auxiliary variable $\beta\in\mathcal B$. The objective of a policy $H$ is
\begin{equation}
\label{eq:costJ}
J(H) = \inf_{\beta \in \mathcal B}\Epath\left[\mathcal L(\beta,X^H,\alpha^H)\right].
\end{equation}
The policy class $\mathcal H$ consists of the measurable feedback maps $H:\Xc\to\Ac$, and the control problem is to choose $H\in\mathcal H$ that minimizes \eqref{eq:costJ}. The loss $\mathcal L$ is specified by the application. Standard expected running and terminal costs are included by taking $\mathcal B$ to be a singleton and $\mathcal L$ to be their sum. For path-dependent problems, this formulation includes deep-hedging risk objectives with an optimized certainty equivalent (OCE) representation \cite{bental2007optimized,buehler2019deep}.

 
\subsection{Random Feature Neural Network Policies}
\label{sec:rfnn_setup}

We now introduce the RFNN policy class with which we approximate the feedback policies of Section~\ref{sec:control_problem}.
An RFNN policy evaluates random features on the policy input $x\in\Xc$ and trains only the output layer.
A single random feature is defined as
\begin{equation}
\label{eq:single_random_feature}
    \phi_{A,b}(x)=\sigma\!\left(A^\top x+b\right),
\end{equation}
where $A\in\R^m$, $b\in\R$, and $\sigma(z)=\max\{z,0\}$ is the ReLU activation used throughout the paper. We sample the hidden weight vectors $(A_k)_{k=1}^N$ and hidden biases $(b_k)_{k=1}^N$ iid within each collection and independently of each other. We write $\Theta_N$ for the sampled hidden parameters $(A_k,b_k)_{k=1}^N$ and $H_N^\Theta$ for the resulting policy.
Two sources of randomness enter the problem, the path randomness generated by the innovations, and the sampled hidden parameters that define the features. To distinguish them, we write $\Epath$ for the path expectation and $\Erf$ for the expectation over the independent RFNN feature draw.
The RFNN policy has $r$ output coordinates, so we write
$H_N^\Theta(x)=\big(H_{N,1}^\Theta(x),\ldots,H_{N,r}^\Theta(x)\big)$. In the shared-feature architecture, all output coordinates use the same collection of random features $\mathcal{R}_N=\{\phi_1,\ldots,\phi_N\}$:
\begin{equation}
\label{eq:shared_feature_policy}
    H_{N,j}^\Theta(x)
    =
    b^{\rm out}_j+\sum_{k=1}^{N} w_{k,j}\phi_k(x),
    \qquad j=1,\ldots,r .
\end{equation}
Here $b^{\rm out}_j$ and $w_{k,j}$ are trainable output parameters. This architecture is equivalent to a single hidden layer neural network in which only the last layer is trainable.
In the coordinate-specific architecture, each output coordinate has its own sampled feature set $\mathcal{R}_{N,j}=\{\phi_{1,j},\ldots,\phi_{N,j}\}$:
\begin{equation}
\label{eq:coordinate_feature_policy}
    H_{N,j}^\Theta(x)
    =
    b^{\rm out}_j+\sum_{k=1}^{N} w_{k,j}\phi_{k,j}(x),
    \qquad j=1,\ldots,r .
\end{equation}
The trainable output parameters $b^{\rm out}_j$ and $w_{k,j}$ have the same role as in the shared-feature case. In both architectures the readout is projected onto $\Ac$ by coordinatewise clipping. The two architectures differ only in whether those features are shared across output coordinates and the appropriate choice depends on the task.

Both architectures have the same trainable output parameters. In the rest of the paper, we do not require independence between features generated for each coordinate, so all results apply to both structures.
\subsection[Assumptions]{Assumptions}
\label{sec:assumptions}

Before moving to the main results, we summarize the assumptions needed for the analysis. We begin with the following assumptions in the approximation framework of \cite{Gonon2023random}.

\begin{assumption}
\label{ass:bounded_policy_input}
The policy input space satisfies
$\Xc\subseteq[-M,M]^m$.
The domain of the RFNN policy is bounded, including both the initial state and every state visited by the admissible coupled trajectories.
\end{assumption}

\begin{assumption}
\label{ass:fourier_prerequisites}
Let $M>0$ and let $H=(H_1,\ldots,H_r):\Xc\to\Ac$. Assume that $H$ admits an extension to $[-M,M]^m$ and that, writing $H$ also for this extension, for each $j=1,\ldots,r$ there exists $G_j:\R^m\to\mathbb C$ such that
\begin{equation}
    H_j(x)=\int_{\R^m} e^{i x\cdot \xi}G_j(\xi)d\xi,
    \qquad x\in[-M,M]^m .
\end{equation}
Suppose that, for each $j$, $\int_{\R^m}\max(1,\|\xi\|^2)|G_j(\xi)|d\xi<\infty$.
Let the hidden weights and biases have strictly positive densities $\pi_A$ on $\R^m$ and $\pi_b$ on $\R$, drawn independently of each other, and define $\bar F(s):=2\int_{-s}^{0}\pi_b(u)^{-1}du$.
Assume that $\bar F(s)\in\R$ for all $s\in\R$ and that, for each $j=1,\ldots,r$, the following constants are finite:
\begin{equation}
\label{eq:scalar_integrability_constant}
\begin{aligned}
I^{\rm int}_j :={}& \max(16,M^2)\int_{\R^m}
\Bigl[\bar F(M\|\xi\|_1)\|\xi\|_1^2 \\
&+\bigl(\bar F(1)-\bar F(-1)\bigr)\max(1,\|\xi\|^2)\Bigr]
\frac{(|G_j(\xi)|+|G_j(-\xi)|)^2}{\pi_A(\xi)}\,d\xi <\infty ,
\end{aligned}
\end{equation}
\begin{equation}
\label{eq:coefficient_bound_constant}
C_{\mathrm w,j}
:=
\sup_{(u,\xi)\in\R\times\R^m}
\bigl(\mathbf 1_{[-M\|\xi\|_1,0]}(u)+4\,\mathbf 1_{[-1,1]}(u)\bigr)
\frac{|G_j(\xi)|+|G_j(-\xi)|}{\pi_b(u)\pi_A(\xi)}
<\infty.
\end{equation}
\end{assumption}

\begin{remark}
Although policies are evaluated only at the decision dates, Assumption~\ref{ass:fourier_prerequisites} is stated for a function of both state and time. It is natural to extend the policy smoothly in the state variable from $\mathcal X$ to the cube $[-M,M]^{m-1}$. Under conditions on the feature distributions, if at every decision date the policy satisfies the assumption as a function of the state on $[-M,M]^{m-1}$, then the extension required in Assumption~\ref{ass:fourier_prerequisites} exists. Full detail is provided in the supplement.
\end{remark}

In addition to the above assumptions, we impose output-parameter bounds and a first-moment condition on the hidden-parameter distribution.

\begin{assumption}
\label{ass:policy_class}
Fix a benchmark policy $H^\star$ satisfying Assumption~\ref{ass:fourier_prerequisites}, and let $C_{\mathrm w,j}$ be its finite constants in \eqref{eq:coefficient_bound_constant}. Fix the feature draw $\Theta_N$ and write $H_\theta$ for the policy with output parameters $\theta=(W,b^{\rm out})$. Write $B_\star:=(\sum_{j=1}^r C_{\mathrm w,j}^2)^{1/2}$, and let $B_{\rm out}$ be a finite positive deterministic radius independent of $N$ and of the feature draw. The output parameters lie in
\begin{equation}
\label{eq:policy_class}
\mathcal W_N=\Big\{(W,b^{\rm out}):\|W_{j,\cdot}\|_2\le\frac{C_{\mathrm w,j}}{\sqrt N}\ \text{for each }j,\ \|b^{\rm out}\|_2\le B_{\rm out}\Big\},
\end{equation}
where $W_{j,\cdot}$ is the $j$th row of $W$. The associated RFNN policy class is $\mathcal H_N:=\{H_\theta:\theta\in\mathcal W_N\}$.
\end{assumption}

\begin{assumption}
\label{ass:feature_law}
The hidden parameter draws $(A_k,b_k)_{k=1}^N$ are iid and have a finite first moment,
\begin{equation}
\label{eq:feature_first_moment}
\Erf\bigl[\|A_1\|_2+|b_1|\bigr]<\infty .
\end{equation}

\end{assumption}
\begin{remark}
In the numerical experiments, the hidden weights and biases are drawn from centered Student-$t$ densities with unit scale and five and two degrees of freedom, respectively, as in \cite{Gonon2023random}. These distributions have strictly positive densities and finite first moments, satisfying the corresponding sampling-density and moment requirements.
\end{remark}

Finally, we have the following assumptions on the controlled transition and the trajectory loss.
\begin{assumption}
\label{ass:stable_coupled_trajectories}
There are deterministic constants $L_X,L_A<\infty$, independent of $t$, $y$, $\eta$ and the policy, such that
\begin{equation}
\label{eq:controlled_update_lipschitz}
\|\Psi_t(y,z,a,\eta)-\Psi_t(y,\tilde z,\tilde a,\eta)\|_1
\le L_X\|z-\tilde z\|_1+L_A\|a-\tilde a\|_1
\end{equation}
for all $t=1,\ldots,T$, $y\in\mathcal Y$, $z,\tilde z\in\mathcal Z$, $a,\tilde a\in\Ac$ and $\eta\in\Ec$.
\end{assumption}

\begin{assumption}
\label{ass:trajectory_loss}
The auxiliary set $\mathcal B\subseteq\R$ is nonempty and compact. The induced pathwise losses are measurable. There are deterministic constants $K_\ell<\infty$ and $0<U<\infty$, independent of the feature width and draw, such that
\begin{equation}
\label{eq:trajectory_loss_lipschitz}
\begin{aligned}
&|\mathcal L(\beta,X,\alpha)-\mathcal L(\widetilde\beta,\widetilde X,\widetilde\alpha)|\\
&\quad\le K_\ell\left(
|\beta-\widetilde\beta|
+\sum_{t=1}^{T+1}\|X_t-\widetilde X_t\|_1
+\sum_{t=1}^{T}\|\alpha_t-\widetilde\alpha_t\|_1
\right),
\end{aligned}
\end{equation}
and
$|\mathcal L(\beta,X,\alpha)|\le U$
for all $\beta,\widetilde\beta\in\mathcal B$, $X_t,\widetilde X_t\in\mathcal X$ and $\alpha_t,\widetilde\alpha_t\in\Ac$.
\end{assumption}
In particular, the population objectives and their empirical counterparts are finite and lie in $[-U,U]$.

\section[Theory]{Theory}
\label{sec:theory}

\subsection[Properties of RFNN]{Properties of RFNN}
\label{sec:rfnn_policies}
Random feature approximations have been studied as computationally efficient alternatives to fully trained neural networks.
We first quote the following uniform approximation theorem from \cite[Theorem~1]{Gonon2023random}.
\begin{theorem}
\label{thm:scalar_rfnn_approx}
Under Assumption~\ref{ass:fourier_prerequisites}, for each $j=1,\ldots,r$, there exist output weights $w_{1,j},\ldots,w_{N,j}$, measurable with respect to the sampled hidden parameters, such that the random feature network
$H_{w,j}(x)=\sum_{k=1}^{N}w_{k,j}\sigma(A_k^\top x+b_k)$
satisfies
\begin{equation}
\label{eq:scalar_rfnn_rate}
\Erf\left[
\sup_{x\in[-M,M]^m}|H_{w,j}(x)-H_j(x)|
\right]
\leq
\frac{4(M\sqrt{m}+1)\sqrt{I^{\rm int}_j}}{\sqrt{N}}.
\end{equation}
Moreover, these weights satisfy, almost surely,
\begin{equation}
\label{eq:scalar_weight_bound}
|w_{k,j}| \leq \frac{C_{\mathrm w,j}}{N},
\qquad k=1,\ldots,N.
\end{equation}
\end{theorem}
Applying the bound of Theorem~\ref{thm:scalar_rfnn_approx} to each coordinate and aggregating the coordinate errors gives the following corollary.
\begin{corollary}
\label{cor:vector_rfnn_approx}
Let $H^\star=(H^\star_1,\ldots,H^\star_r):\Xc\to\Ac$ be a benchmark policy satisfying Assumption~\ref{ass:fourier_prerequisites}. For each coordinate $j$, let $I^{\rm int}_j$ be the integrability constant in \eqref{eq:scalar_integrability_constant} associated with $H^\star_j$ and its Fourier transform, and let $C_j=4(M\sqrt{m}+1)\sqrt{I^{\rm int}_j}$ be the constant of Theorem~\ref{thm:scalar_rfnn_approx} applied to coordinate $j$. Then there exists an RFNN approximation $H_N^\Theta$ with measurable output weights and zero output bias such that
\begin{equation}
\label{eq:vector_rfnn_rate}
\Erf\left[
\sup_{x\in\Xc}\|H_N^\Theta(x)-H^\star(x)\|_1
\right]
\leq
\frac{rC_{\max}}{\sqrt{N}},
\qquad
C_{\max}=\max_{j=1,\ldots,r} C_j .
\end{equation}
The output weight matrix $W_N^\star\in\R^{r\times N}$ of the same approximation satisfies
$\|W_{N,j}^\star\|_2\le C_{\mathrm w,j}N^{-1/2}$ for each $j=1,\ldots,r$, which is the row-level form of \eqref{eq:scalar_weight_bound},
and hence
\begin{equation}
\label{eq:output_weight_frobenius}
\|W_N^\star\|_F\le\frac{B_\star}{\sqrt N},
\qquad
B_\star=\Big(\sum_{j=1}^{r}C_{\mathrm w,j}^2\Big)^{1/2}<\infty,
\end{equation}
where $C_{\mathrm w,j}$ is the constant in \eqref{eq:coefficient_bound_constant}.
\end{corollary}
Assumption~\ref{ass:policy_class} places $H_N^\Theta$ in $\mathcal H_N$ almost surely, which is the reason why we restrict the analysis to the bounded ball \eqref{eq:policy_class}.  We then prove that every policy in this class is Lipschitz in its input and in its output parameters. For the later proofs, we further show that the logarithms of the resulting Lipschitz constants grow at most logarithmically in $N$ in expectation over the feature draw.
\begin{proposition}
\label{prop:class_lipschitz}
Let $\Theta_N$ be a fixed feature draw and let $H_\theta\in\mathcal H_N$. For all policy inputs $x,x'\in\Xc$,
\begin{equation}
\label{eq:class_state_lipschitz}
\|H_\theta(x)-H_\theta(x')\|_1\le L_{\mathrm{class}}\|x-x'\|_1,
\qquad
L_{\mathrm{class}}:=\,B_\star\Bigl(\frac1N\sum_{j=1}^r\sum_{k=1}^N\|A_{k,j}\|_2^2\Bigr)^{1/2}.
\end{equation}

\end{proposition}

\begin{proposition}
\label{prop:class_discrepancy}
Under Assumption~\ref{ass:bounded_policy_input}, let $\Theta_N$ be a fixed feature draw. For all $\theta,\theta'\in\mathcal W_N$,
\begin{equation}
\label{eq:class_discrepancy}
\sup_{x\in\Xc}\|H_\theta(x)-H_{\theta'}(x)\|_1
\le
\Lambda_W\,\|W-W'\|_F+\Lambda_{\rm out}\,\|b^{\rm out}-b'^{\rm out}\|_2 ,
\end{equation}
where
\begin{equation}
\label{eq:class_discrepancy_constants}
\Lambda_W:=\Bigl(\sum_{j=1}^r\sum_{k=1}^N\sup_{x\in\Xc}|\phi_{k,j}(x)|^2\Bigr)^{1/2},
\qquad
\Lambda_{\rm out}:=\sqrt r .
\end{equation}
\end{proposition}

\begin{proposition}
\label{prop:constant_integrability}
Under Assumptions~\ref{ass:bounded_policy_input} and~\ref{ass:feature_law}, the weight constant of
\eqref{eq:class_discrepancy_constants} and the state-Lipschitz constant of
\eqref{eq:class_state_lipschitz} satisfy
\begin{equation}
\label{eq:constant_log_growth}
\Erf\bigl[\log\bigl(1+\Lambda_W\bigr)\bigr]\le C_0+\log N,
\qquad
\Erf\bigl[\log\bigl(1+L_{\mathrm{class}}\bigr)\bigr]\le C_0+\log N,
\end{equation}
where $C_0$ depends only on $m,M,r,B_\star$ and the law of the hidden parameters.
\end{proposition}

\subsection[Approximation error]{Approximation error}
\label{sec:approximation_theory}
In this subsection, we convert the policy approximation bound into a bound on the control objective. A policy perturbation changes both the coupled state and action trajectories entering the terminal loss. We therefore propagate the policy approximation error through the controlled dynamics before applying a final loss stability bound. Throughout this subsection, the RFNN feature draw is fixed unless an $\Erf$ expectation is written explicitly.
We first consider two policies $H_1$ and $H_2$ which operate in the same environment, so their trajectories differ only through their own feedback laws.
\begin{theorem}
\label{thm:policy_error_propagation}
Under Assumptions~\ref{ass:bounded_policy_input} and~\ref{ass:stable_coupled_trajectories}, let $H_1$ and $H_2$ be any two policies on the same innovations. Suppose that $H_1$ is Lipschitz in the state, uniformly in the decision date, with $\|H_1(X,t)-H_1(\tilde X,t)\|_1\le L_1\|X-\tilde X\|_1$ for all states $X,\tilde X$ and all dates $t$, and that the two policies satisfy the uniform discrepancy bound $\varepsilon:=\sup_{x\in\Xc}\|H_1(x)-H_2(x)\|_1<\infty$. Set $\Gamma=L_X+L_A L_1$.
Then, on every common realization of the innovations, for each $t=1,\ldots,T+1$,
\begin{equation}
\label{eq:state_error_propagation}
    \|X_t^{H_1}-X_t^{H_2}\|_1
    \leq
    L_A\varepsilon\sum_{s=0}^{t-2}\Gamma^s,
\end{equation}
and for $t=1,\ldots,T$,
\begin{equation}
\label{eq:action_error_propagation}
    \|\alpha_t^{H_1}-\alpha_t^{H_2}\|_1
    \leq
    \varepsilon
    \left(
    1+L_1 L_A\sum_{s=0}^{t-2}\Gamma^s
    \right).
\end{equation}
Empty sums are zero.
\end{theorem}
The next proposition turns these state and action bounds into a bound on the objectives.
\begin{proposition}
\label{prop:final_loss_stability}
Under Assumption~\ref{ass:trajectory_loss}, the policy objective satisfies the following joint state and action stability bound along a common path:
\begin{equation}
\label{eq:final_loss_stability}
|J(H_1)-J(H_2)|\le K_\ell\Bigl[\sum_{t=1}^{T+1}\Epath\|X_t^{H_1}-X_t^{H_2}\|_1+\sum_{t=1}^{T}\Epath\|\alpha_t^{H_1}-\alpha_t^{H_2}\|_1\Bigr].
\end{equation}

\end{proposition}
Overall, we connect the policy error with the change in the loss.
\begin{corollary}
\label{cor:class_objective_stability}
Under Assumptions~\ref{ass:bounded_policy_input}, \ref{ass:stable_coupled_trajectories} and \ref{ass:trajectory_loss}, let $H_1$ be $L_1$-Lipschitz in the state, uniformly in time, and let $H_2$ be a second policy with finite uniform discrepancy $\varepsilon=\sup_{x\in\Xc}\|H_1(x)-H_2(x)\|_1$. Then
\begin{equation}
\label{eq:class_objective_stability}
|J(H_1)-J(H_2)|\le K_T(L_1)\varepsilon.
\end{equation}
Here $\Gamma=L_X+L_A L_1$ and
\begin{equation}
\label{eq:combined_cost_constant}
K_T(L_1):=K_\ell\left[T+L_A(1+L_1)\sum_{s=0}^{T-1}(T-s)\Gamma^s\right].
\end{equation}
The same propagation gives the pathwise estimate
\begin{equation}
\label{eq:combined_path_cost}
\sup_{\beta\in\mathcal B}
|\mathcal L(\beta,X^{H_1},\alpha^{H_1})-\mathcal L(\beta,X^{H_2},\alpha^{H_2})|
\le K_T(L_1)\varepsilon.
\end{equation}
\end{corollary}
Corollary~\ref{cor:class_objective_stability} converts a uniform policy discrepancy into a discrepancy of objectives. We apply it with $H_1=H^\star$ and $H_2=H_N^\Theta$. The Fourier condition on $H^\star$ implies that $H^\star$ is Lipschitz, so the resulting estimate combines with Corollary~\ref{cor:vector_rfnn_approx} to give the objective approximation bound below.
\begin{theorem}
\label{thm:rfnn_objective_approx}

Suppose that Assumptions~\ref{ass:bounded_policy_input}, \ref{ass:policy_class}, \ref{ass:stable_coupled_trajectories} and \ref{ass:trajectory_loss} hold. Let $H^\star=(H^\star_1,\ldots,H^\star_r):\Xc\to\Ac$ be a benchmark policy satisfying Assumption~\ref{ass:fourier_prerequisites}, and let $H_N^\Theta$ be the RFNN approximation supplied by Corollary~\ref{cor:vector_rfnn_approx}. Let $L_\star$ be a Lipschitz constant of $H^\star$ in the state, uniformly in the decision date.  The RFNN policy then satisfies the objective approximation bound
\begin{equation}
\label{eq:rfnn_loss_error_bound}
    \Erf\left[
    |J(H_N^\Theta)-J(H^\star)|
    \right]
    \leq
    K_T(L_\star)
    \frac{rC_{\max}}{\sqrt{N}}.
\end{equation}
\end{theorem}
\begin{remark}
\label{rem:stability_constant_growth}
For a fixed horizon, $K_T$ grows at most polynomially in the problem data. With $\Gamma=L_X+L_AL$ the geometric sum is at most $T^2\max(1,\Gamma)^{T-1}$, hence
\begin{equation}
\label{eq:stability_constant_polynomial}
K_T(L)\le K_\ell\bigl[T+L_A(1+L)T^2\max(1,\Gamma)^{T-1}\bigr],
\end{equation}
which, up to the factor $K_\ell$, is a polynomial of total degree at most $2T$ in $L_A$, $L_X$ and $L$, and of degree at most $T$ in $L$ when $L_A$ and $L_X$ are fixed. The stability step preserves the polynomial dimension dependence of the approximation rate provided $L_\star$ and the constants $I^{\mathrm{int}}_j$ grow polynomially in the dimension.
\end{remark}
 

\subsection[Statistical error]{Statistical error}
\label{sec:statistical_error}
We now study the statistical error arising from replacing the population objective by its empirical counterpart based on finitely many simulated paths. Let $\omega_1,\ldots,\omega_n$ be i.i.d. draws of the innovation sequence defined in Section~\ref{sec:control_problem}, independent of the RFNN feature draw. The empirical objective is
\begin{equation}
\label{eq:empirical_objective}
\widehat J_n(H)
= \inf_{\beta\in\mathcal B}
\frac1n\sum_{i=1}^n\mathcal L\bigl(\beta,X^H(\omega_i),\alpha^H(\omega_i)\bigr),
\end{equation}
and the corresponding constrained empirical risk minimizer is
$\widehat H\in\arg\min_{H\in\mathcal H_N}\widehat J_n(H)$.

Define the pathwise loss $\ell(\omega;\theta,\beta):=\mathcal L\bigl(\beta,X^{H_\theta}(\omega),\alpha^{H_\theta}(\omega)\bigr)$.
We first prove that it is Lipschitz in the parameters.

\begin{proposition}
\label{prop:param_loss_lipschitz}
Suppose that Assumptions~\ref{ass:bounded_policy_input}, \ref{ass:stable_coupled_trajectories} and \ref{ass:trajectory_loss} hold. Then, for every realization $\omega$ of the innovations, every $\beta\in\mathcal B$, and all $\theta,\theta'\in\mathcal W_N$,
\begin{equation}
\label{eq:param_loss_lipschitz}
\begin{aligned}
|\ell(\omega;\theta,\beta)-\ell(\omega;\theta',\beta)|
&\le
L_\ell^W\|W-W'\|_F+L_\ell^{\rm out}\|b^{\rm out}-b'^{\rm out}\|_2,\\
L_\ell^W&:=\Lambda_W K_T(L_{\rm class}),
\qquad
L_\ell^{\rm out}:=\Lambda_{\rm out} K_T(L_{\rm class}),
\end{aligned}
\end{equation}
with $K_T$ defined in \eqref{eq:combined_cost_constant}.
\end{proposition}

The preceding Lipschitz estimate gives a covering-number bound for the pathwise loss class. An entropy integral bounds its Rademacher complexity \cite[Theorem~5.22]{Wainwright2019}, and symmetrization as in the proof of \cite[Theorem~8]{BartlettMendelson2002} yields the following uniform empirical deviation bound.

\begin{theorem}
\label{thm:estimation_error}
Suppose that Assumptions~\ref{ass:bounded_policy_input}, \ref{ass:policy_class}, \ref{ass:stable_coupled_trajectories} and \ref{ass:trajectory_loss} hold. Fix the feature draw $\Theta_N$ and let $\Epath$ denote expectation over the i.i.d. path draws conditional on $\Theta_N$. Set $r_N:=r(N+1)$. There is a constant $C_{\mathrm{est}}$, independent of $n$, $N$ and the feature draw, such that
\begin{equation}
\label{eq:estimation_error_bound}
\Epath\Bigl[\sup_{H\in\mathcal H_N}\bigl|\widehat J_n(H)-J(H)\bigr|\Bigr]
\le
\frac{\kappa_{\mathrm{est}}}{\sqrt n}
+
\frac{2U}{\sqrt n},
\qquad
\kappa_{\mathrm{est}}
:=
C_{\mathrm{est}}U\sqrt{r_N+1}\,
\sqrt{\log\Bigl(2+\frac{\mathcal L_\ell}{U}\Bigr)}.
\end{equation}
where $\mathcal L_\ell:=L_\ell^W+L_\ell^{\rm out}=(\Lambda_W+\Lambda_{\rm out})K_T(L_{\rm class})$. The coefficient $\kappa_{\mathrm{est}}$ depends on the fixed feature draw through $\mathcal L_\ell$.
\end{theorem}

By Proposition~\ref{prop:constant_integrability}, we can bound the expectation of the logarithmic term in $\kappa_{\mathrm{est}}$ and hence the expected growth of $\kappa_{\mathrm{est}}$ with the feature width.

\begin{proposition}
\label{prop:estimation_constant_growth}
Under the assumptions of Theorem~\ref{thm:estimation_error} and Assumption~\ref{ass:feature_law}, there are constants $C_{\log}$ and $C_{\rm grow}$, independent of $n$ and $N$, such that
\begin{equation}
\label{eq:loss_log_bound}
\Erf\Bigl[\log\Bigl(2+\frac{\mathcal L_\ell}{U}\Bigr)\Bigr]
\le
C_{\log}\bigl(1+\log N\bigr),
\end{equation}
and consequently
\begin{equation}
\label{eq:kappa_growth}
\Erf[\kappa_{\mathrm{est}}]
\le
C_{\rm grow}\sqrt{N\bigl(1+\log N\bigr)}
=
O\bigl(\sqrt{N\log N}\bigr).
\end{equation}
\end{proposition}

We finally combine the approximation and statistical errors. Let $H_N^\star$ be a minimizer of $J$ over $\mathcal H_N$. For an exact empirical minimizer $\widehat H$ over $\mathcal H_N$, the overall error satisfies
\begin{equation}
\label{eq:excess_risk_decomposition}
J(\widehat H)-J(H^\star)
\le
J(H_N^\star)-J(H^\star)
+
\bigl|\widehat J_n(\widehat H)-J(\widehat H)\bigr|
+
\bigl|\widehat J_n(H_N^\star)-J(H_N^\star)\bigr|.
\end{equation}
The first term is the approximation error controlled by Theorem~\ref{thm:rfnn_objective_approx}. The latter two terms are the empirical--population deviations at the two policies, and each is bounded above by the statistical error controlled by Theorem~\ref{thm:estimation_error}. Together, these bounds yield the following corollary.
\begin{corollary}
\label{cor:excess_risk}
Under the hypotheses of Theorems~\ref{thm:rfnn_objective_approx} and~\ref{thm:estimation_error}, and Assumption~\ref{ass:feature_law}, suppose that the draws $\omega_1,\ldots,\omega_n$ are iid and independent of the feature bank. Assume that a measurable exact empirical minimizer over $\mathcal H_N$ exists, and let $\widehat H$ be such a minimizer. Then
\begin{equation}
\label{eq:excess_risk_bound}
\Erf\,\Epath\bigl[J(\widehat H)-J(H^\star)\bigr]
\le
C\left(
\frac{1}{\sqrt N}
+
\sqrt{\frac{N(1+\log N)}{n}}
\right),
\end{equation}
where $C$ is independent of $N$ and $n$.
\end{corollary}

The two terms in \eqref{eq:excess_risk_bound} represent the approximation and estimation errors, respectively. Increasing the number of random features therefore improves the approximation error while increasing the estimation error. Balancing the two terms gives $N^2(1+\log N)\asymp n$ and hence, asymptotically, $N\asymp\sqrt{n/\log n}$.
Substituting this choice into \eqref{eq:excess_risk_bound} yields
\begin{equation}
\Erf\,\Epath\bigl[J(\widehat H)-J(H^\star)\bigr]
\le
C\left(\frac{\log n}{n}\right)^{1/4}.
\end{equation}

\section{Numerical experiments}
\label{sec:numerical_error_decomposition}

In this section, we study how feature width, sample size and feature selection affect the performance of RFNN policies. We first introduce the Markovian control and Heston hedging examples. After describing the training and evaluation procedures, we check the assumptions of the analysis on a bounded formulation of each example. We then examine approximation and statistical errors and the resulting trade-off in feature width, and compare RFNN with fully trained networks. We also evaluate Parsimonious Random Feature Learning (PRFL), which replaces the fixed RFNN feature set with a dictionary that repeatedly grows, is trained and is pruned. We compare its performance with RFNN and with fully trained networks at matched retained width and training budget. Throughout, $N$ denotes feature width, $n$ is the number of paths in the empirical objective, and lower objective values indicate better performance.

\subsection{Benchmark problems}
\label{sec:controlled_experiments}
We consider the Markovian control problems studied in \cite{bachouch2022deep}, including semilinear and linear--quadratic control. We also consider hedging under the Heston model, following the setup of \cite{buehler2019deep}.

\subsubsection{Markovian control}

\paragraph{Semilinear control in dimension $d=100$.}
\label{sec:num_semilinear_example}

The first Markovian example is associated with the semilinear Hamilton--Jacobi--Bellman equation
\begin{equation}
\label{eq:num_semilinear_hjb}
\partial_t v(t,X)+\Delta v(t,X)-|Dv(t,X)|^2=0,
\qquad
v(1,X)=g_{\mathrm{PDE}}(X).
\end{equation}
Using $-|p|^2=\inf_{a\in\R^d}\{|a|^2+2a\cdot p\}$, the corresponding Euler-discretized controlled state satisfies
\begin{equation}
\label{eq:num_semilinear_dynamics}
X_{t+1}=X_t+2h\alpha_t+\sqrt{2h}\,\varepsilon_{t+1},
\end{equation}
where the innovations are independent standard Gaussian vectors. The policy minimizes the expected pathwise loss
\begin{equation}
\label{eq:num_semilinear_objective}
J_{\mathrm{SL}}(H)
=\Epath\left[h\sum_{t=1}^{20}\|\alpha_t^H\|_2^2
+g_{\mathrm{PDE}}(X_{21}^H)\right].
\end{equation}
We take $d=100$, $X_1=0$, $h=1/20$, and $g_{\mathrm{PDE}}(X)=\log\bigl(\tfrac12(1+|X|^2)\bigr)$. Following \cite{bachouch2022deep}, the Hopf--Cole transformation gives the value function
\begin{equation}
\label{eq:num_hopf_cole_reference}
v(t,x)
=
-\log\E\left[
\exp\bigl(-g_{\mathrm{PDE}}(x+\sqrt{2}W^{\rm SL}_{1-t})\bigr)
\right]
\end{equation}
for $(t,x)\in[0,1]\times\R^d$, and the benchmark policy is the optimal feedback of the continuous-time problem, $\alpha_t=-Dv(t,X_t)$. 

\paragraph{Linear--quadratic control.}
\label{sec:num_lq_example}

The second Markovian example has state dynamics
\begin{equation}
\label{eq:num_lq_dynamics}
X_{t+1}
=
X_t+h(BX_t+\Gamma_{\rm LQ}\alpha_t)
+\sum_{j=1}^{p}D_j\alpha_t\sqrt h\,\varepsilon_{t+1}^j
\end{equation}
and objective
\begin{equation}
\label{eq:num_lq_objective}
J_{\mathrm{LQ}}(H)
=\Epath\left[h\sum_{t=1}^{20}\left((X_t^H)^\top QX_t^H+\lambda_{\mathrm{LQ}}\|\alpha_t^H\|_2^2\right)
+(X_{21}^H)^\top PX_{21}^H\right].
\end{equation}
The benchmark uses $Q=P=I_d$, $B=I_d$, $\Gamma_{\rm LQ}=\mathbf 1_d$, $D_j=e_j$, $p=d$, $h=1/20$, and $X_1=0.1\mathbf 1_d$, with $d\in\{10,100\}$ and $\lambda_{\mathrm{LQ}}=1$. Here $e_j$ is the $j$th coordinate vector in $\mathbb R^d$. The innovations $\varepsilon_{t+1}=(\varepsilon_{t+1}^1,\ldots,\varepsilon_{t+1}^p)$ are independent standard Gaussian.

Following \cite{bachouch2022deep}, the value function remains quadratic for the discrete simulator, $V_t(X)=X^\top \Pi_tX$, and the matrices are computed backwards from $\Pi_{21}=P$. Writing $\mathbf C=h\Gamma_{\rm LQ}$, the recursion is
\begin{equation}
\label{eq:num_lq_riccati}
\Pi_t
=
hQ+(I+hB)^\top
\left(
\Pi_{t+1}
-
\frac{\Pi_{t+1}\mathbf C\mathbf C^\top \Pi_{t+1}}
{\lambda_{\mathrm{LQ}}h+\mathbf C^\top \Pi_{t+1}\mathbf C
+h\operatorname{tr}(\Pi_{t+1})}
\right)
(I+hB).
\end{equation}
The LQ benchmark policy is the linear feedback $\alpha_t=-K_tX_t$, where $\Pi_t$ follows the Riccati recursion \eqref{eq:num_lq_riccati} and $K_t=(\lambda_{\mathrm{LQ}}h+\mathbf C^\top\Pi_{t+1}\mathbf C+h\operatorname{tr}(\Pi_{t+1}))^{-1}\mathbf C^\top\Pi_{t+1}(I+hB)$. 

\subsubsection{Heston hedging}
\label{sec:num_nocost_heston_example}

The financial example uses a Heston model for the underlying asset $S^1$ and its stochastic variance $v$:
\begin{equation}
\label{eq:num_heston_model}
\begin{aligned}
dS^1_t &= \mu_{\rm H}S^1_t\,dt+\sqrt{v_t}S^1_t\,dW^{S}_t,\\
dv_t &= \kappa_{\rm H}(\bar v_{\rm H}-v_t)\,dt+\sigma_{{\rm H},v}\sqrt{v_t}\,dW^{v}_t,\\
d\langle W^{S},W^{v}\rangle_t &= \varrho\,dt.
\end{aligned}
\end{equation}
Time is measured in years, with step $h=1/365$ and $T=30$ decision dates, so the horizon is $Th=30/365$. Decisions are made at dates $1,\ldots,T$, the state is observed at the $T+1$ dates $0,h,\ldots,Th$, and decision date $t$ corresponds to calendar time $(t-1)h$. The market coordinates of the state are $(S^1_t,v_t)$. The variance is not directly tradable, so a variance swap is included as the second hedging instrument. Its price, together with the underlying, is a functional of the state trajectory,
\begin{equation}
\label{eq:num_tradable_map}
\mathcal S_t(X_{[1:t]})
=
\left(
S^1_t,
\int_0^{(t-1)h} v_s\,ds
+\frac{v_t-\bar v_{\rm H}}{\kappa_{\rm H}}\left(1-e^{-\kappa_{\rm H}(T+1-t)h}\right)
+\bar v_{\rm H}(T+1-t)h
\right),
\end{equation}
where the integrated variance is approximated by the trapezoidal rule. The strategy chooses holdings $\delta_t\in\Ac\subseteq\R^2$ in the underlying and the variance swap. For a European call payoff $\mathrm{Payoff}=(S^1_{T+1}-K)_+$, its terminal P\&L is
\begin{equation}
\label{eq:num_heston_pnl}
\mathrm{PnL}_{\lambda_{\mathrm{tc}}}(X,\boldsymbol\delta)
=\sum_{t=1}^{T}\delta_t^\top(\mathcal S_{t+1}-\mathcal S_t)-\mathrm{Payoff}
-\lambda_{\mathrm{tc}}\sum_{t=1}^{T-1}\mathcal S_{t+1}^\top|\delta_{t+1}-\delta_t|.
\end{equation}
The simulations use $\mu_{\rm H}=0$, $S^1_0=100$, $v_0=0.04$, $\kappa_{\rm H}=1$, $\bar v_{\rm H}=0.04$, $\sigma_{{\rm H},v}=2$, $\varrho=-0.7$, $K=100$, and CVaR level $q=0.5$.

Following \cite{buehler2019deep}, the policies are evaluated under square loss and CVaR, each with an auxiliary offset $\beta$:
\begin{align}
J_{\mathrm{sq}}(H)
&=\min_{\beta\in\R}\Epath\left[\left(\beta-\mathrm{PnL}_{\lambda_{\mathrm{tc}}}(X^H,\boldsymbol\delta^H)\right)^2\right],\label{eq:num_square_loss}\\
J_{\mathrm{CVaR},q}(H)
&=\inf_{\beta\in\R}\Epath\left[\beta+\frac{1}{1-q}\left(-\mathrm{PnL}_{\lambda_{\mathrm{tc}}}(X^H,\boldsymbol\delta^H)-\beta\right)_+\right].\label{eq:num_cvar_loss}
\end{align}

The offset $\beta$ is a trainable parameter during training, and it is re-optimized on the test data during evaluation. The policy state is $X_t=(S^1_t,v_t,\delta_{t-1})\in\R^4$, with feedback holdings $\delta_t=H(X_t,t)$. The network input replaces $S^1_t$ by $\log S^1_t$, and the previous holding is omitted when $\lambda_{\mathrm{tc}}=0$. The market coordinates follow the simulated path, which is common to all policies, and the holding is updated by the action. The no-cost case is obtained by setting $\lambda_{\mathrm{tc}}=0$.

\subsection{Training and evaluation}
\label{sec:numerical_training}

Following Section~\ref{sec:rfnn_setup}, RFNN trains the output coefficients of a fixed random-feature dictionary, while FTNN trains both layers. Following the theoretical analysis, let $H^\star$ denote the benchmark of each example in Section~\ref{sec:controlled_experiments}, $H_N^\star$ the policy that minimizes $J$ over $\mathcal H_N$, and $\widehat H_{n,N}$ the policy that minimizes $\widehat J_n$ over $\mathcal H_N$. The benchmark policy is specified in each example. For Heston hedging it is the best observed policy on the evaluated objective across the tested widths, seeds and repetitions.

\paragraph{Adam training.} Both RFNN and FTNN are trained with the Adam optimizer, with a learning rate that decays during training. The training schemes differ by problem. In the Markovian examples, minibatches contain $1024$ paths for LQ and $10^4$ for semilinear control, or the entire sample when $n$ is smaller, and training runs for $200$ epochs for LQ and $50$ epochs for semilinear control, with $50$ updates per epoch. The initial learning rate $10^{-3}$ is held for $30$ epochs and then multiplied by $0.25$ every five epochs, with these milestones four times longer for LQ. In Heston hedging, minibatches contain up to $10^4$ paths, training runs for $1600$ epochs with $10$ updates per epoch, and the learning rate $0.005$ is multiplied by $0.1$ every $400$ epochs. To stabilize training, we apply weight decay $10^{-3}$ in the Markovian examples. In Heston hedging, we apply gradient clipping at $1$, and the inputs are normalized and bounded. For the reference policies in the Markovian examples, a new sample of paths is generated at every update. Training becomes less stable as the width grows in both RFNN and FTNN, and in the semilinear example a few clearly deviating runs are omitted from the reported averages.   

\paragraph{Direct regression with fixed random features.} With no transaction costs and square loss, the previous holding can be omitted from the policy state. Fixing the hidden features makes the terminal P\&L affine in the output coefficients and allows the RFNN policy to be fitted directly by constrained least squares. Collecting the simulated trading gains and payoff into a design matrix $\mathbf A$ and response vector $Y$ gives
\begin{equation}
\label{eq:num_nocost_regression}
\widehat J_n(\theta)
=
\frac1n\|Y-\mathbf A\theta\|_2^2.
\end{equation}
We solve this problem over the prescribed bounded output-layer domain, adjusting the constraint multiplier to enforce feasibility. The reference policy $H_N^\star$ is computed by the same constrained regression on an independent sample of $10^6$ paths. 

\paragraph{Evaluation and repetition protocol.} The approximation and statistical-error experiments use three feature seeds for RFNN or three initialization seeds for FTNN. For finite-sample experiments, the policy is retrained on $3$ to $20$ training samples for each seed, with more repetitions at smaller sample sizes. The PRFL comparison uses ten training seeds for the Markovian examples and three for Heston. Each final objective is estimated by Monte Carlo on $10^6$ independent evaluation paths with the policy held fixed.

\subsection{Assumptions in the numerical examples}
\label{sec:num_assumption_scope}

The Gaussian innovations and Heston price and variance processes considered above have unbounded support, whereas our convergence analysis is established on bounded policy-input domains. To check Assumption~\ref{ass:bounded_policy_input}, we consider bounded formulations of these examples, obtained by truncating the driving innovations or observed market paths at fixed levels and restricting controls to a compact admissible action set. Over the finite horizon, the resulting state dynamics yield bounded reachable states. In the numerical experiments, we do not actually do the truncation or the projection, as we can choose the bounds that are large enough during training.

Assumption~\ref{ass:fourier_prerequisites} requires a density condition on the feature draws and a regularity condition on the benchmark policies. The hidden weights and biases are centered unit-scale Student-$t$ draws with five and two degrees of freedom, respectively, whose densities are strictly positive, so the density condition holds. The benchmark policy itself need not satisfy this regularity condition. On the bounded formulation, any Lipschitz benchmark is within a small distance of a smooth policy, which we may keep $\Ac$-valued by enlarging the action bound. Corollary~\ref{cor:class_objective_stability} makes the objective change at most a constant multiple of that distance, so policies satisfying Assumption~\ref{ass:fourier_prerequisites} attain almost the same performance and the reported gaps change by a small amount. Feature draws are independent with finite first moment, so Assumption~\ref{ass:feature_law} holds as well.

Assumption~\ref{ass:policy_class} requires the output parameters to lie within prescribed bounds independent of the network width. The Heston regression experiments enforce these bounds through constrained regression, while the remaining experiments use Adam with weight decay. The bounds are estimated from the reference policy and decay at the rate $N^{-1/2}$.

In each example, the expression inside the expectation defining the objective is the pathwise loss. The Markovian examples use no auxiliary optimization, so $\mathcal B=\{0\}$. For Heston hedging, the bounded P\&L allows the auxiliary offset $\beta$ to be restricted to a compact interval. On the bounded domains considered above, the transitions and pathwise losses are Lipschitz, which yields finite constants for Assumptions~\ref{ass:stable_coupled_trajectories} and~\ref{ass:trajectory_loss}. For Heston hedging the controlled block is the holding, which the action replaces, so the assumption holds with $L_X=0$ and $L_A=1$ on the unnormalized state, and the market map need not be Lipschitz.

\subsection{RFNN approximation and statistical errors}

For an exact empirical minimizer $\widehat H_{n,N}$, the decomposition in \eqref{eq:excess_risk_decomposition} writes the objective gap as the approximation term and the empirical--population deviations at the two policies,

\begin{equation}
\label{eq:num_excess_decomposition}
J(\widehat H_{n,N})-J(H^\star)
\le
J(H_N^\star)-J(H^\star)
+
D^{\mathrm{emp}}_{n,N}
+
D^{\mathrm{ind}}_{n,N},
\end{equation}
where
$D^{\mathrm{emp}}_{n,N}=\bigl|\widehat J_n(\widehat H_{n,N})-J(\widehat H_{n,N})\bigr|$ and $D^{\mathrm{ind}}_{n,N}=\bigl|\widehat J_n(H_N^\star)-J(H_N^\star)\bigr|$.
The first term is the loss from restricting the policy to $\mathcal H_N$, and it is bounded by a constant multiple of $N^{-1/2}$. Each of the last two terms is at most the uniform empirical--population deviation over $\mathcal H_N$, which Theorem~\ref{thm:estimation_error} bounds in expectation by a constant multiple of $\sqrt{N(1+\log N)/n}$. 

\begin{figure}[!htbp]
\centering
\includegraphics[width=0.84\textwidth]{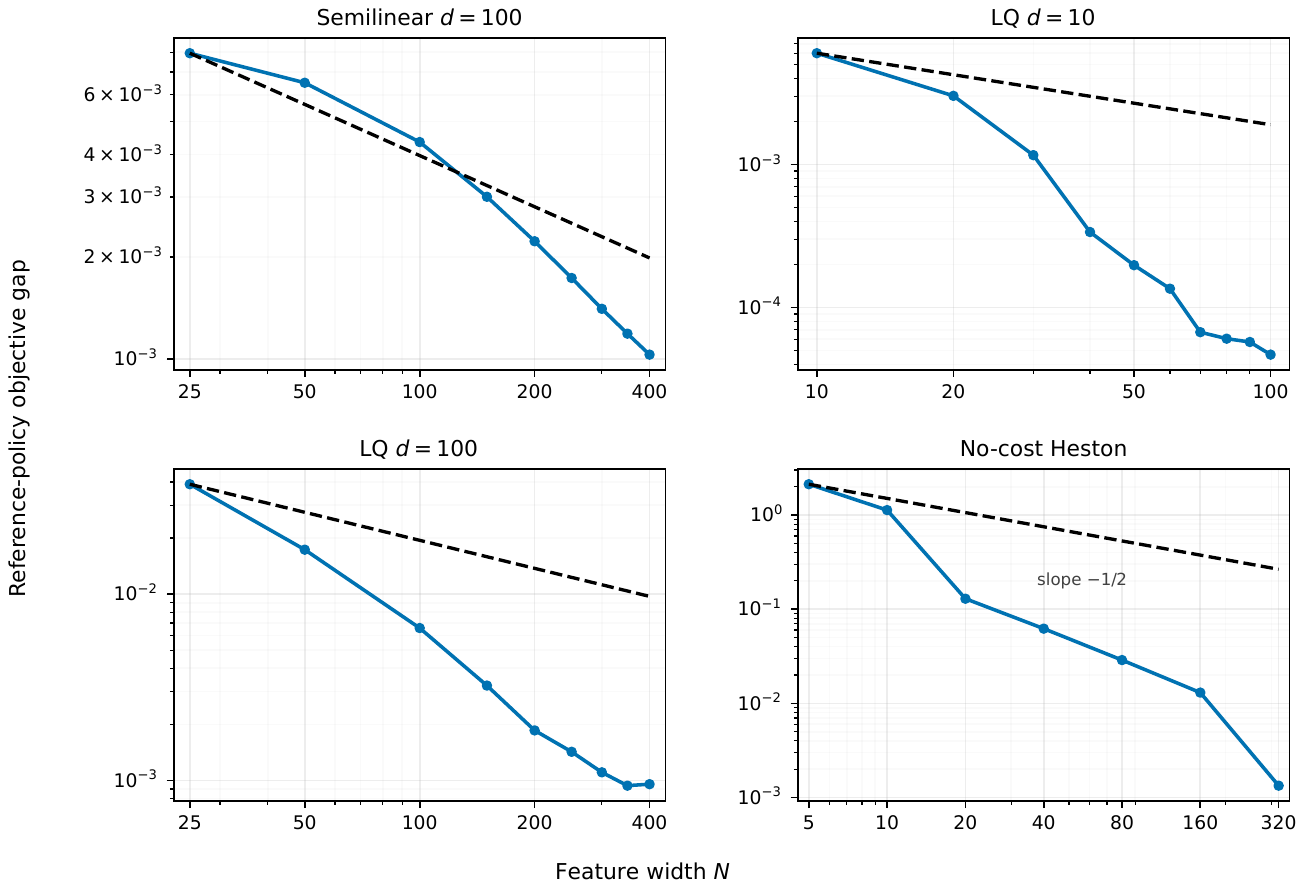}
\caption{Reference-policy objective gaps versus feature width, where the dashed guide marks the $N^{-1/2}$ rate.}
\label{fig:approximation_error_controlled}
\end{figure}

\begin{figure}[!htbp]
\centering
\includegraphics[width=0.84\textwidth]{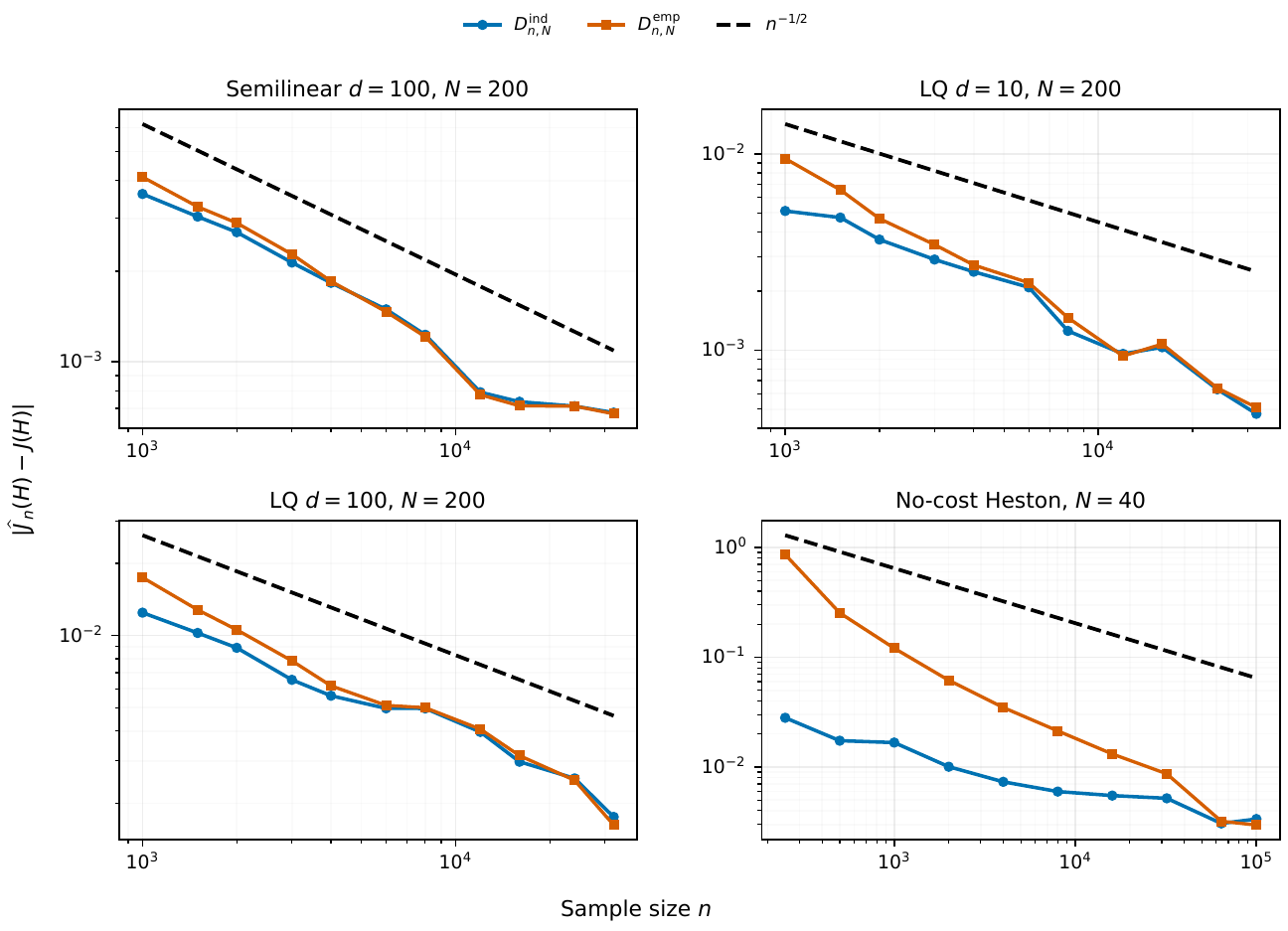}
\caption{Sample-size dependence of the empirical--population deviation at the sample-independent policy $H_N^\star$ and the empirical policy $\widehat H_{n,N}$.}
\label{fig:statistical_sample_size}
\end{figure}

We first measure each term on its own. Figure~\ref{fig:approximation_error_controlled} shows the approximation term. The plotted gap $J(H_N^\star)-J(H^\star)$ compares the reference policy with the benchmark. The gap decreases over the width grid in all four problems, and the $N^{-1/2}$ guide marks the rate of Theorem~\ref{thm:rfnn_objective_approx}. The measured gaps fall faster than the guide in all four problems. On the other hand, Figure~\ref{fig:statistical_sample_size} shows the statistical term at fixed width. Both deviations decrease with the sample size in all four panels. Since both policies lie in $\mathcal H_N$, each deviation is bounded above by the uniform error of Theorem~\ref{thm:estimation_error}, which decays like $n^{-1/2}$ and matches the guide in the figure.  Overall, both terms match the theoretical results.

\begin{figure}[!htbp]
\centering
\includegraphics[width=\textwidth]{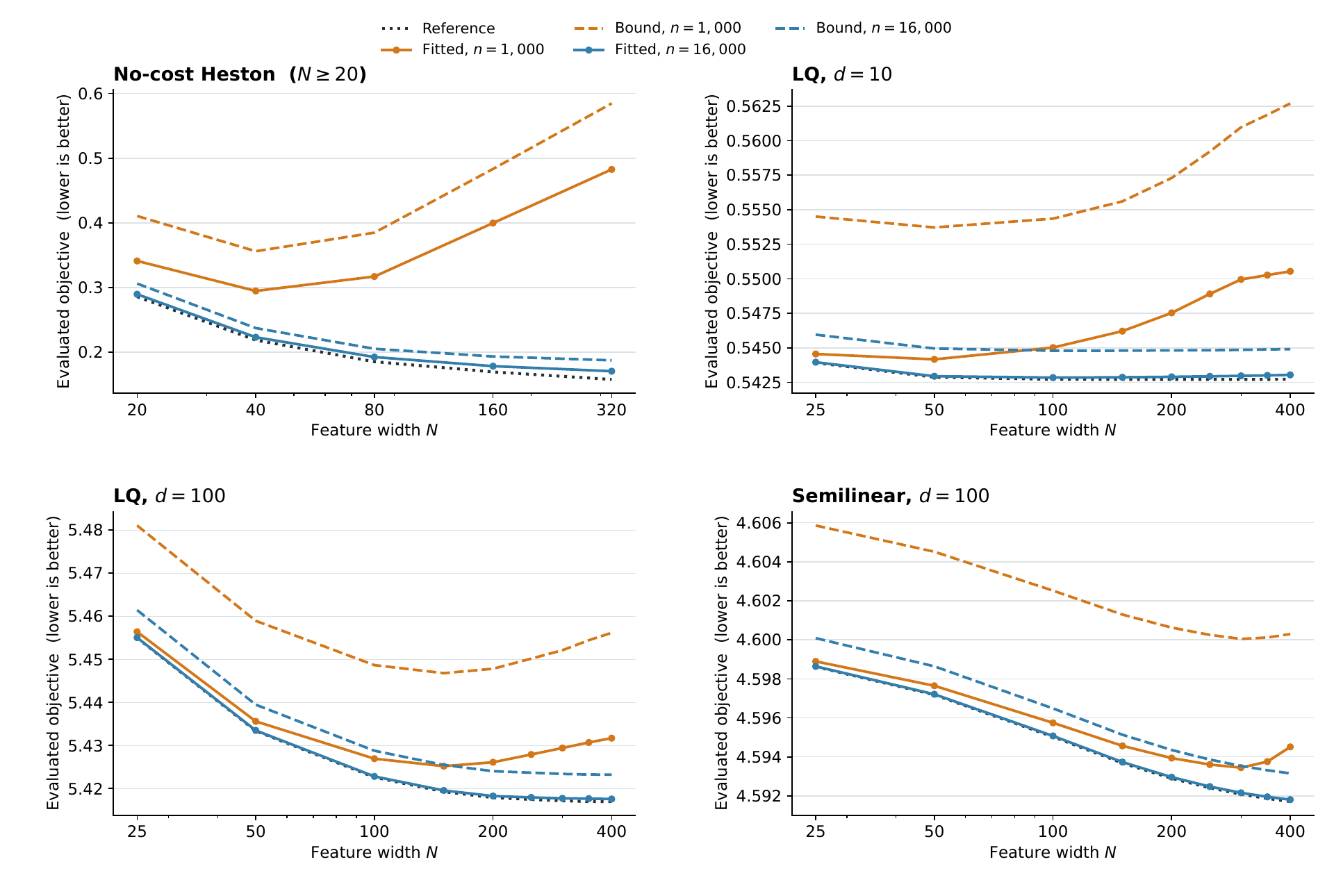}
\caption{Reference and fitted objectives versus feature width for $n=1000$ and $n=16000$. The solid curves are the fitted objective and the dashed curves the right-hand side of \eqref{eq:num_excess_decomposition}, both up to the constant $J(H^\star)$. The Heston panel displays $N\geq20$.}
\label{fig:width_tradeoff_comparison}
\end{figure}

Figure~\ref{fig:width_tradeoff_comparison} brings the two errors together as the width increases. The dashed curves trace the right-hand side of \eqref{eq:num_excess_decomposition} and the solid curves the left-hand side, both up to the constant $J(H^\star)$. Because the reference objective decreases with $N$ while the two deviations grow, the trade-off already appears in the dashed curves. The solid curves follow the same trend and remain close to them, with the gap narrowing as the sample size grows. At $n=1000$ the fitted objective reaches a minimum before turning upward over part of the grid, and at $n=16000$ it continues to decrease, so the width with the lowest fitted objective moves to larger $N$. The figure shows $n=1000$ and $n=16000$ as two examples of this balance.

\subsection{RFNN versus fully trained networks}
\label{sec:numerical_results}
\label{sec:rfnn_ftnn_sample_size}

This subsection compares RFNN with a fully trained neural network (FTNN), which trains both layers, on the same training samples, feature widths, sample-size grids and evaluation protocol. The comparison shows how the two policy classes behave as the number of training paths changes. 

\begin{figure}[!htbp]
\centering
\includegraphics[width=0.84\textwidth]{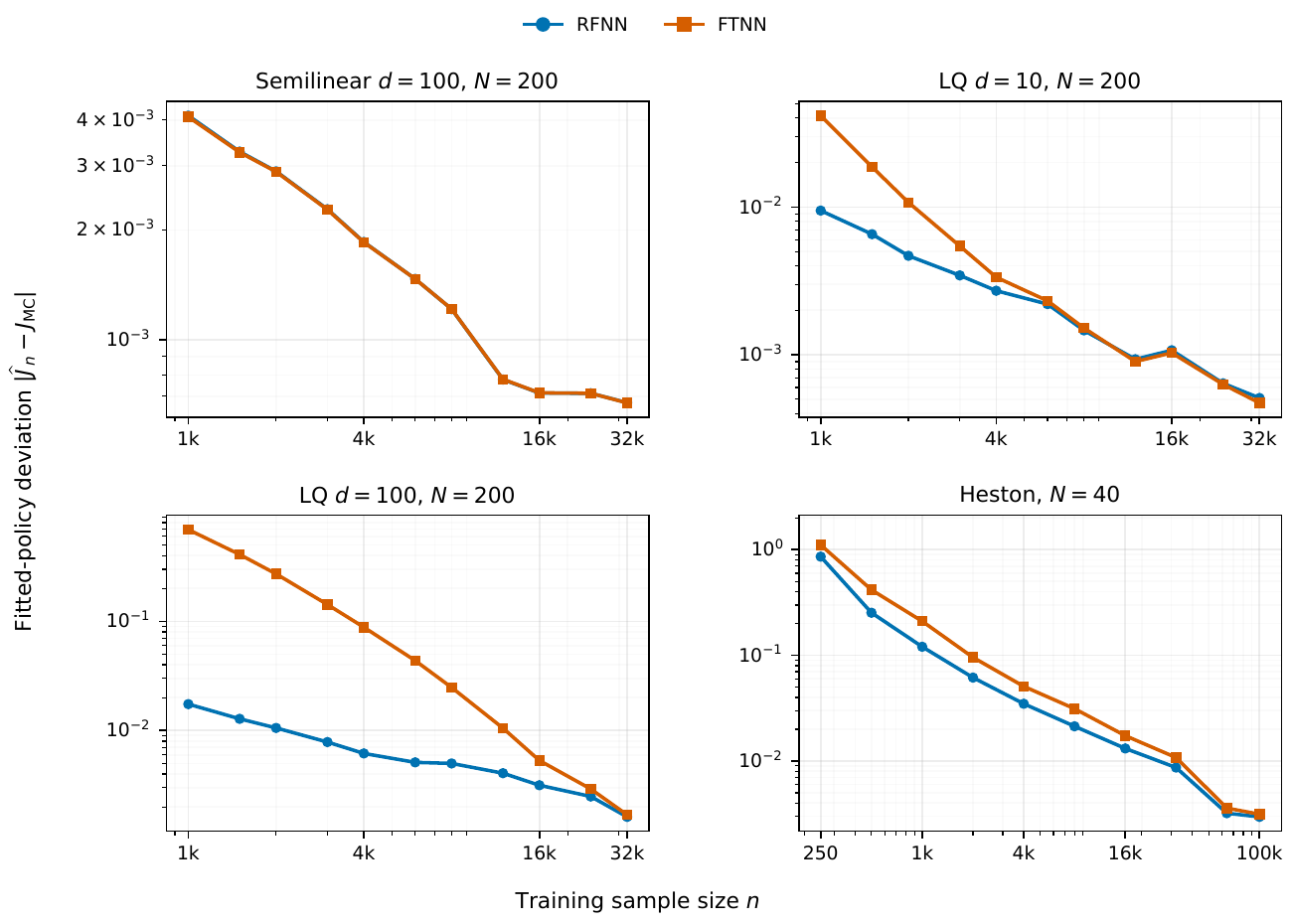}
\caption{Empirical--population deviations of the empirical policy for RFNN and FTNN versus training sample size.}
\label{fig:rfnn_ftnn_statistical_sample_size}
\end{figure}

With limited training data, RFNN can benefit from the smaller policy class obtained by fitting only the output coefficients of a fixed feature representation. Figure~\ref{fig:rfnn_ftnn_statistical_sample_size} compares the statistical error and illustrates this small-sample advantage in the LQ and Heston examples. The deviations generally decrease and the difference between methods narrows as the sample size grows. The semilinear results are already similar across the tested range. The zero control attains $4.6002$ and the benchmark $4.5907$, so the total improvement available in this example is below $10^{-2}$, which leaves little room for the two methods to differ.

\begin{figure}[!htbp]
\centering
\includegraphics[width=0.84\textwidth]{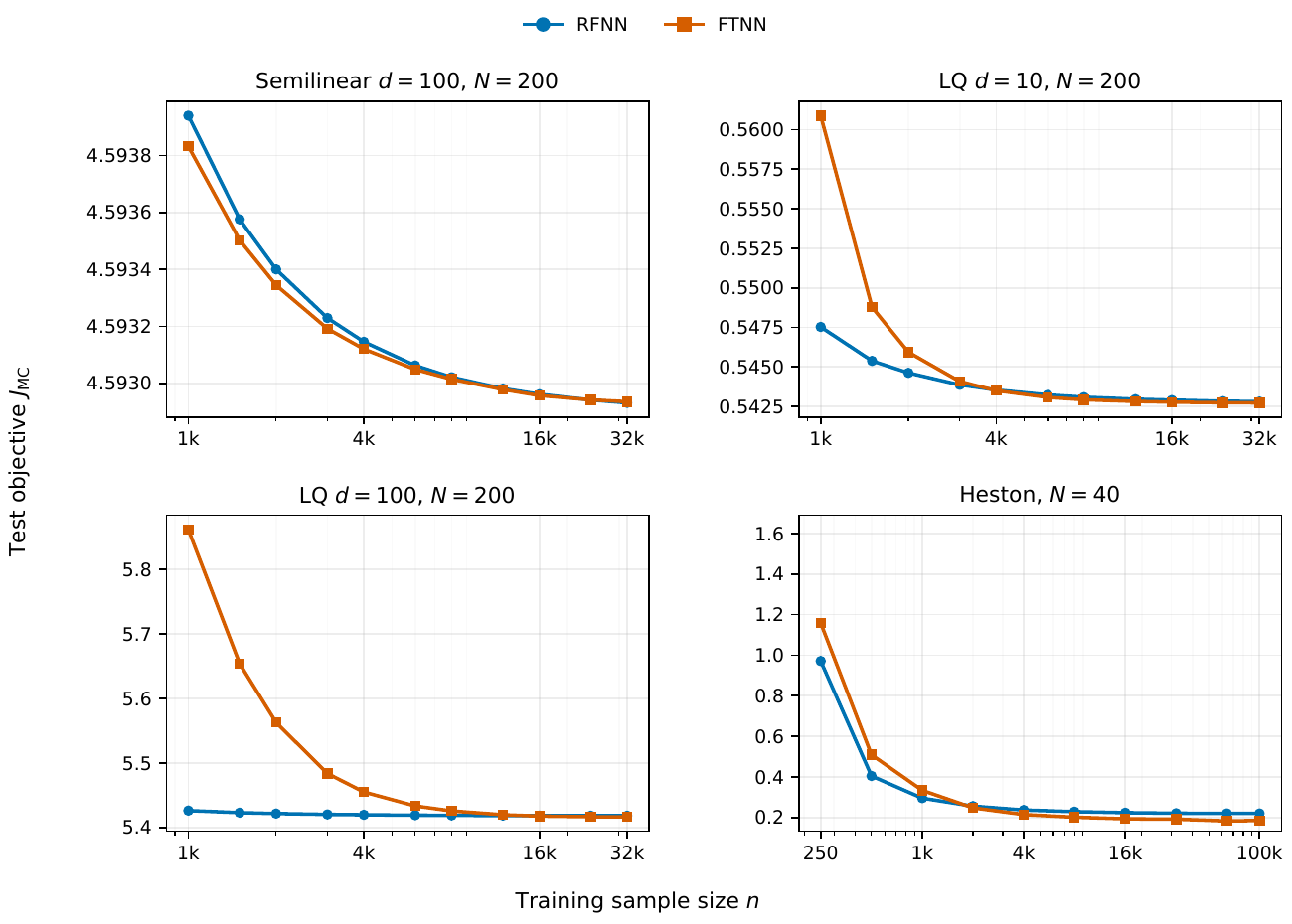}
\caption{Independent test objectives of RFNN and FTNN versus training sample size.}
\label{fig:rfnn_ftnn_test_objective_sample_size}
\end{figure}

Figure~\ref{fig:rfnn_ftnn_test_objective_sample_size} gives the corresponding test objectives, and the same mechanism explains their ordering. With few training paths the statistical error is relatively large, so the smaller RFNN class is an advantage, and this advantage carries through to the test objective. As the sample grows the statistical error falls, the approximation advantage of the larger FTNN class becomes more important, and FTNN overtakes. The crossover is visible in the two LQ examples and in Heston. In semilinear control the two policies are already close at every tested sample size, so the crossover does not appear. Overall, the advantage of RFNN over FTNN is limited to small samples, and FTNN matches or beats it once enough paths are available.

\subsection{Adaptive feature selection}
\label{sec:adaptive_feature_experiments}

Parsimonious Random Feature Learning (PRFL) builds on the online learning scheme of \cite{lemoal2025parsimonious}, which was developed for online regularized regression and used inside backward dynamic programming. Here we consider it under a different framework, direct policy optimization. The idea we keep is that the feature set itself changes during training. We test whether this improves policy performance at a given retained architecture while training only the output coefficients.

\subsubsection{PRFL method and pruning rule}
\label{sec:prfl}

The feature set changes during training in $R$ cycles. PRFL starts from an initial dictionary of width $s_{\rm add}$. Each cycle trains the output layer, scores the current features and deletes the features with the smallest weights allowed by the pruning tolerance. The first $R-1$ cycles then add a fresh block of $s_{\rm add}$ features, and after the last cycle training continues on the retained dictionary. The procedure therefore inspects $R s_{\rm add}$ features, while only part of them is retained in the final policy.

The pruning rule controls how much one pruning step can change the trained output layer. For the shared dictionaries used in the Markovian benchmarks, write the trained policy on an enlarged dictionary $E$ of $N$ features as
$H_E(x)=b^{\rm out}+\sum_{\phi_k\in E}w_k\phi_k(x)$,
where $b^{\rm out}\in\R^r$ is the bias and $w_k=(w_{k,1},\ldots,w_{k,r})$ is the vector of output weights for feature $\phi_k$.  PRFL scores each feature by the squared weights across the output coordinates, $s_k=\sum_{j=1}^{r}|w_{k,j}|^2$, and deletes the longest prefix of the increasing-score order whose cumulative score stays within the pruning tolerance. With $D$ the deleted set and $W=(w_1,\ldots,w_N)$, the pruning budget is exactly
$\|W-W_{\mathrm{pruned}}\|_F^2=\sum_{k\in D}s_k\leq\eta_{\mathrm{pr}}$.
This coefficient-space sparsification parallels the random-feature projection used in POLRF~\cite{lemoal2025parsimonious}, where orthonormal feature coordinates likewise make the squared projection error the sum of the squared deleted coefficients.
The smallest scores have the smallest sum, so this rule deletes as many features as the budget allows. Zeroing columns preserves the parameter bound and leaves the projection unchanged, hence $H_{\mathrm{pruned}}\in\mathcal H_N$ whenever $H_E\in\mathcal H_N$. By Proposition~\ref{prop:class_discrepancy}, deleting these features changes the policy by at most $\Lambda_W\sqrt{\eta_{\mathrm{pr}}}$ in uniform discrepancy. Corollary~\ref{cor:class_objective_stability} then bounds the change in the objective by
$|J(H_{\mathrm{pruned}})-J(H_E)|\le K_T(L_{\mathrm{class}})\Lambda_W\sqrt{\eta_{\mathrm{pr}}}$.
 The right-hand side grows with the width through its constants, but a small enough pruning tolerance keeps the perturbation from each deletion under control. In the numerical experiment below, we show that this control is enough for the gain from inspecting more features to dominate. Convergence and optimality of the full procedure remain open for future research.

\subsubsection{PRFL experiment}
\label{sec:prfl_experimental_protocol}

In each PRFL run, the adaptive phase occupies the first half of training, and the retained dictionary is trained for the remaining half. The adaptive phase lasts $100$ epochs for LQ, $25$ for semilinear control and $800$ for Heston hedging, with $20$, $5$ and $50$ epochs per cycle, so a run lasts $200$, $50$ and $1600$ epochs in total. Each epoch makes $50$ updates for the Markovian examples and $10$ for Heston, giving $10000$, $2500$ and $16000$ updates per run.  The pruning tolerance varies by problem, at $10^{-5}$ for the Markovian examples, $10^{-3}$ for Heston at zero transaction cost, and $10^{-2}$ at positive costs.

The Markovian controls use initial and added feature counts $s_{\rm add}\in\{5,10,20\}$ for semilinear control, $\{5,10,\ldots,25\}$ for LQ with $d=10$, and $\{5,10,\ldots,40\}$ for LQ with $d=100$. Each epoch uses $50$ newly simulated minibatches of $1024$ paths for all three methods, and test objectives and retained widths are averaged over ten seeds. The Heston comparison covers square loss and CVaR at transaction-cost rates $\lambda_{\mathrm{tc}}\in\{0,0.0005,0.005\}$, with $s_{\rm add}\in\{25,50,75,100\}$ per instrument and three training seeds per configuration, using newly simulated minibatches for all methods. The remaining configuration follows Section~\ref{sec:numerical_training}, and the full PRFL setup is given in the supplement.

PRFL is trained with Adam, and one optimizer is kept across the add-train-prune cycles. When features are deleted, the Adam state of the retained output weights is carried over, including the step count and the first and second moment estimates. When features are added, the new output weights enter with zero moment estimates while the existing ones keep theirs. This avoids resetting the optimizer state of surviving features and the transient change in their updates that a reset would cause.

\subsubsection{Results}
\label{sec:prfl_benchmark_results}

\begin{figure}[!htbp]
\centering
\includegraphics[width=\textwidth]
{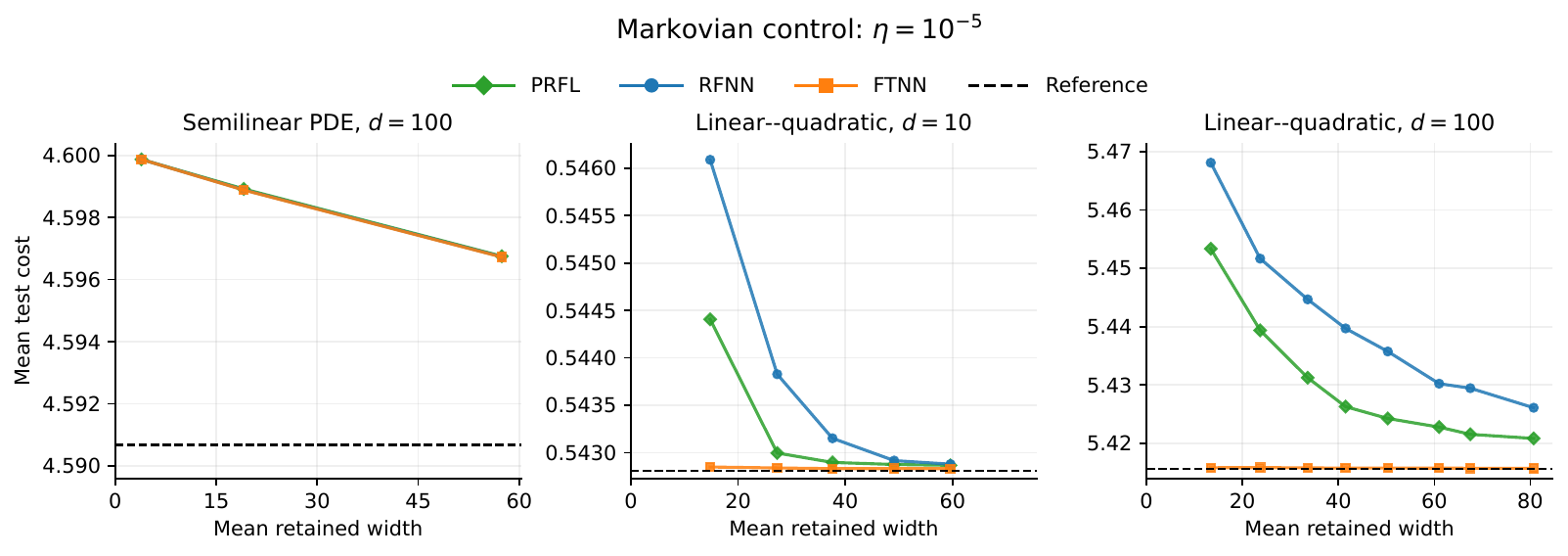}
\caption{Test objectives versus mean retained width for RFNN, FTNN and PRFL in the Markovian benchmarks; dashed lines indicate the corresponding reference values.}
\label{fig:sc_benchmarks}
\end{figure}

\begin{figure}[!htbp]
\centering
\includegraphics[width=\textwidth]
{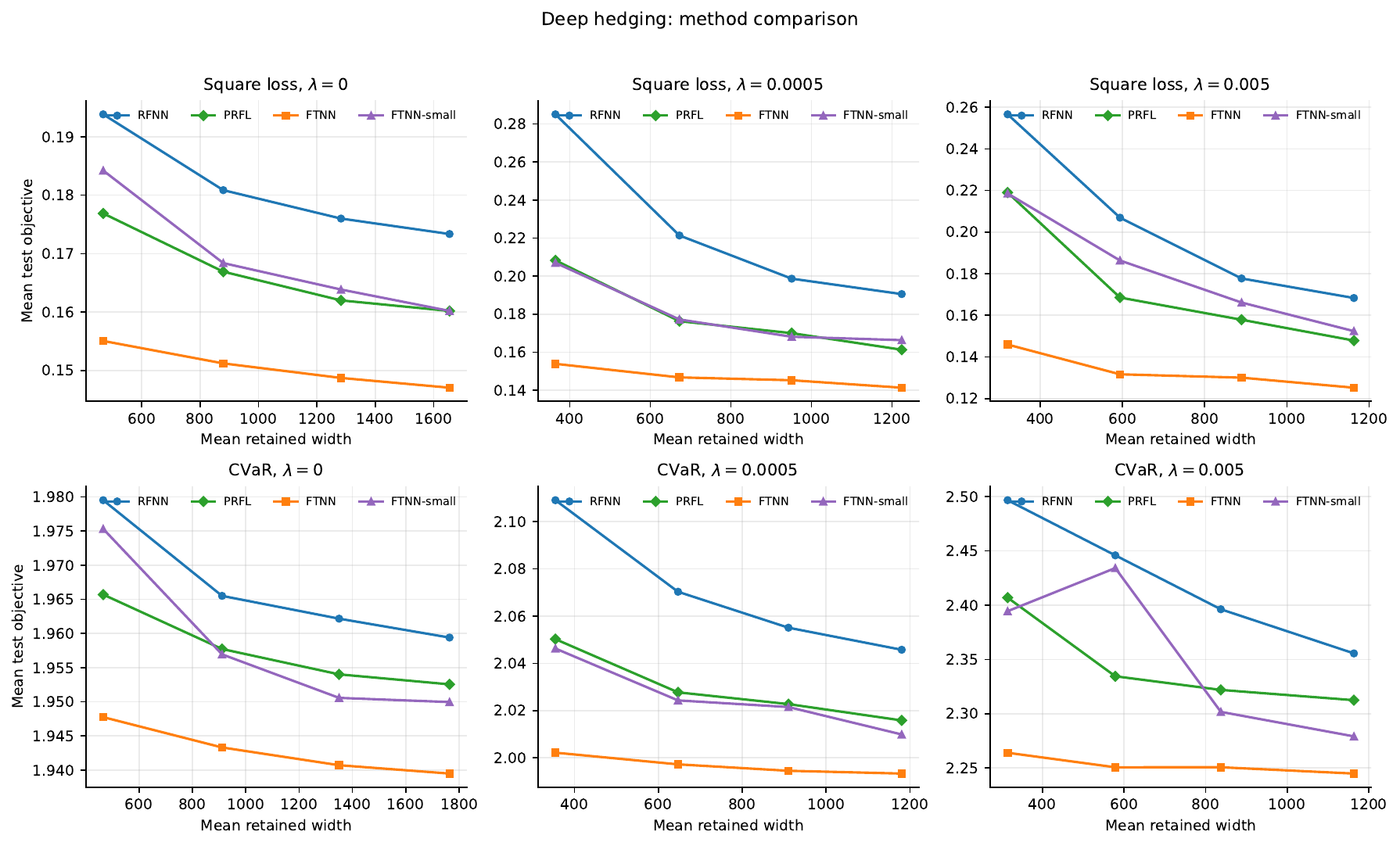}
\caption{Test objectives in Heston hedging under square loss and CVaR at three transaction-cost rates. The horizontal axis is the mean retained width of the PRFL policy, and FTNN-small is plotted at the same position although its width is smaller.}
\label{fig:tc_method_comparison}
\end{figure}

For each PRFL run, we compare PRFL with RFNN and FTNN trained from scratch with the same retained width, network architecture and number of training epochs. In Heston, the widths are matched separately for each hedging instrument. Figure~\ref{fig:sc_benchmarks} shows the Markovian controls and Figure~\ref{fig:tc_method_comparison} the Heston comparison, with the retained width on the horizontal axis and the test objective on the vertical axis. PRFL improves on the matched RFNN in the LQ and Heston comparisons and narrows the gap to the fully trained FTNN. Inspecting more features therefore overcomes both the loss from feature deletion and the loss from interrupting training when features are added and removed. We perform an additional ablation study in the supplement to further support this.  In semilinear control RFNN and FTNN already perform similarly, and PRFL achieves comparable results. PRFL and RFNN have fewer trainable parameters than a fully trained network of the same width, so the Heston experiment also includes FTNN-small, whose width approximately matches their trainable-parameter count. Relative to FTNN-small, PRFL reaches a very similar performance across the figures. Overall, PRFL therefore provides a meaningful extension of RFNN, keeping the output-layer training while matching a fully trained network of comparable trainable size.

\bibliographystyle{siamplain}
\bibliography{reference}

\appendix
\section{Proofs}
\label{app:proofs}
This appendix collects the proofs of the preceding theoretical results.

\subsection{Discussion of Assumption~\ref{ass:fourier_prerequisites}}
\label{app:student_time_extension}
We construct the extension for one output coordinate. Let $f_t$ and $G_t$ be the policy and its Fourier density at date $t$. Write $d=d_X$.

Choose a real $q\in C_c^\infty((-1/2,1/2))$ with $q(0)=1$, and let $\rho$ be its Fourier density. Define
$R(\xi)=\bigl(1+\|\xi\|^2\bigr)^{1/2}$.
Define the extension by
\begin{equation}
\widetilde H(x,u)
=\sum_{t=1}^T\int_{\mathbb R^d}e^{ix\cdot\xi}G_t(\xi)
q\bigl(R(\xi)(u-t)\bigr)\,d\xi.
\end{equation}
Since $R(\xi)\geq1$, the support of $q$ gives $\widetilde H(x,s)=f_s(x)$ on $[-M,M]^d$ for every decision date $s$. Its Fourier density is
\begin{equation}
\label{eq:student_extension_fourier}
\widetilde G(\xi,\tau)
=\sum_{t=1}^T G_t(\xi)e^{-it\tau}
\frac{1}{R(\xi)}
\rho\!\left(\frac{\tau}{R(\xi)}\right).
\end{equation}
By the rapid decay of $\rho$,
\begin{equation}
\begin{aligned}
&\int_{\mathbb R^{d+1}}(1+\|\xi\|^2+|\tau|^2)|\widetilde G(\xi,\tau)|\,d\xi\,d\tau\\
&\quad\leq
\left(\int_{\mathbb R}(1+|z|^2)|\rho(z)|\,dz\right)
\sum_{t=1}^T\int_{\mathbb R^d}(1+\|\xi\|^2)|G_t(\xi)|\,d\xi
<\infty.
\end{aligned}
\end{equation}
where we set $z=\tau/R(\xi)$.   

For the Student-$t$ distributions used in this paper, we verify the remaining two conditions. Write $\pi_d$ and $\pi_{d+1}$ for the densities of $t_\nu(0,I_d)$ and $t_\nu(0,I_{d+1})$, and let $\pi_b$ be the same Student-$t$ bias density in both cases. Set $z=\tau/R(\xi)$ and $a=(\nu+d+1)/2$. Increasing the dimension changes the weight density by
\begin{equation}
\label{eq:student_extension_ratio}
\frac{1}{\pi_{d+1}(\xi,R(\xi)z)}
\leq \frac{C R(\xi)}{\pi_d(\xi)}(1+z^2)^a.
\end{equation}
The factor $R(\xi)$ is why we included $R(\xi)^{-1}$ in $\widetilde G$. For the coefficient condition, the two factors cancel, leaving the original ratio $G_t/\pi_d$ multiplied by $(1+z^2)^a\rho(z)$. For the integral condition, $\widetilde G$ is squared, and the resulting $R(\xi)^{-2}$ cancels with the $R(\xi)$ from the density bound and the $R(\xi)$ from the change of variables $d\tau=R(\xi)\,dz$.

Let the Student-$t$ bias have $\mu$ degrees of freedom, and set $s=\|\xi\|_1$. Since $1/\pi_b(v)\asymp(1+|v|)^{\mu+1}$, the function $\bar F$ satisfies $1+\bar F(u)\asymp(1+u)^{\mu+2}$ for $u\geq0$. Using $\tau=R(\xi)z$ and $R(\xi)\asymp1+s$ gives
\begin{equation}
\label{eq:student_bias_integral_comparison}
1+\bar F\bigl(M(s+|\tau|)\bigr)
\asymp
\bigl(1+\bar F(Ms)\bigr)(1+|z|)^{\mu+2}.
\end{equation}
Similarly, with $D(s):=\sup_{v\in[-Ms,0]\cup[-1,1]}1/\pi_b(v)$,
\begin{equation}
\label{eq:student_bias_coefficient_comparison}
D(s+|\tau|)\asymp D(s)(1+|z|)^{\mu+1}.
\end{equation} Together with the norm factors in the integrability condition, these comparisons bound the joint bias factors by their date-specific counterparts times polynomials in $|z|$. The resulting integrals and suprema are finite because $\rho$ decays faster than any polynomial. Thus, if $I_t$ and $C_{{\rm w},t}$ are the finite constants in Assumption~\ref{ass:fourier_prerequisites} for $f_t$, the corresponding joint constants satisfy
\begin{equation}
\widetilde I\leq C\sum_{t=1}^T I_t,
\qquad
\widetilde C_{\rm w}\leq C\sum_{t=1}^T C_{{\rm w},t}<\infty.
\end{equation}
The condition $\bar F(s)\in\mathbb R$ is unchanged. Hence $\widetilde H$ satisfies Assumption~\ref{ass:fourier_prerequisites} for these Student-$t$ distributions.

The same construction applies to other distributions whenever the upper bounds in \eqref{eq:student_extension_ratio}, \eqref{eq:student_bias_integral_comparison} and \eqref{eq:student_bias_coefficient_comparison} hold with polynomial factors in $|z|$. Under these conditions, the cancellations above and the rapid decay of $\rho$ give the remaining two bounds.

\subsection{Proof of Corollary~\ref{cor:vector_rfnn_approx}}
\label{app:proof:vector_rfnn_approx}
\begin{proof}
By Assumption~\ref{ass:fourier_prerequisites}, Theorem~\ref{thm:scalar_rfnn_approx} applies to each coordinate $H^\star_j$ of the benchmark policy. Write $\mathrm{raw}_j$ for the unclipped readout of coordinate $j$, so that $H^\Theta_{N,j}=\mathrm{clip}(\mathrm{raw}_j)$, where $\mathrm{clip}$ is the coordinatewise projection onto $\Ac$. For every $j$ the coordinate-wise approximation error obeys
\begin{equation}
\label{eq:coordinate_rfnn_rate}
    \Erf\left[
    \sup_{x\in[-M,M]^m}|H^\star_j(x)-\mathrm{raw}_j(x)|
    \right]
    \leq
    \frac{C_j}{\sqrt{N}}.
\end{equation}
On the policy input space the benchmark policy is $\Ac$-valued, so $\mathrm{clip}(H^\star_j(x))=H^\star_j(x)$ for every $x\in\Xc$. Since clipping is nonexpansive,
\begin{equation}
\label{eq:clipping_step}
    |H^\Theta_{N,j}(x)-H^\star_j(x)|
    =
    |\mathrm{clip}(\mathrm{raw}_j(x))-\mathrm{clip}(H^\star_j(x))|
    \leq
    |\mathrm{raw}_j(x)-H^\star_j(x)|
\end{equation}
for every $x\in\Xc$. Suprema over $\Xc$ and expectations therefore inherit \eqref{eq:coordinate_rfnn_rate}. Since the $\ell^1$ norm splits over the coordinates,
\begin{equation}
\label{eq:vector_policy_l1_bound}
\sup_{x\in\Xc}\|H_N^\Theta(x)-H^\star(x)\|_1
\leq
\sum_{j=1}^r \sup_{x\in\Xc} |H_{N,j}^\Theta(x)-H^\star_j(x)| ,
\end{equation}
and taking $\Erf$ expectations of both sides yields \eqref{eq:vector_rfnn_rate}.  For the output weights, \eqref{eq:scalar_weight_bound} applied to coordinate $j$ gives $|w_{k,j}|\le C_{\mathrm w,j}/N$ for $k=1,\ldots,N$, hence $\sum_{k=1}^N w_{k,j}^2\le C_{\mathrm w,j}^2/N$ and
\begin{equation}
\|W_N^\star\|_F^2
=
\sum_{j=1}^r\sum_{k=1}^N w_{k,j}^2
\leq
\frac1N\sum_{j=1}^r C_{\mathrm w,j}^2 .
\end{equation}
\end{proof}

\subsection{Proof of Proposition~\ref{prop:class_lipschitz}}
\label{app:proof:class_lipschitz}
\begin{proof}
Since the projection is nonexpansive, it suffices to bound the output before projection. The one-Lipschitz property of ReLU gives $|\phi_{k,j}(x)-\phi_{k,j}(x')|\le\|A_{k,j}\|_2\|x-x'\|_2$. The output bias cancels. Cauchy--Schwarz first over the features in each coordinate and then over the output coordinates gives
\begin{equation}
\begin{aligned}
&\|H_\theta(x)-H_\theta(x')\|_1\le\|W\|_F\Big(\sum_{j=1}^r\sum_{k=1}^N\|A_{k,j}\|_2^2\Big)^{1/2}\|x-x'\|_2.
\end{aligned}
\end{equation}
Since  $\|x-x'\|_2\le\|x-x'\|_1$ and $\sum_{j=1}^r\|W_{j,\cdot}\|_2^2\le\sum_{j=1}^rC_{\mathrm w,j}^2/N=B_\star^2/N$, the bound \eqref{eq:class_state_lipschitz} follows.
\end{proof}

\subsection{Proof of Proposition~\ref{prop:class_discrepancy}}
\label{app:proof:class_discrepancy}
\begin{proof}

By Assumption~\ref{ass:bounded_policy_input}, $\sup_{x\in\Xc}|\phi_{k,j}(x)|\le M\|A_{k,j}\|_1+|b_{k,j}|<\infty$, so $\Lambda_W$ is finite. Since the projection is nonexpansive, it suffices to bound the output before projection. Write the policy difference in coordinate $j$ as $(b^{\rm out}_j-b'^{\rm out}_j)+\sum_{k=1}^N(w_{k,j}-w'_{k,j})\phi_{k,j}(x)$. Cauchy--Schwarz within each row and then across rows gives
\begin{equation}
\sum_{j=1}^r\left|\sum_{k=1}^N(w_{k,j}-w'_{k,j})\phi_{k,j}(x)\right|
\le
\begin{aligned}
&\Bigl(\sum_{j=1}^r\sum_{k=1}^N\sup_{x\in\Xc}|\phi_{k,j}(x)|^2\Bigr)^{1/2} \|W-W'\|_F .
\end{aligned}
\end{equation}
The coefficient is $\Lambda_W$ in \eqref{eq:class_discrepancy_constants}. For the bias term, $\|b^{\rm out}-b'^{\rm out}\|_1\le\sqrt r\|b^{\rm out}-b'^{\rm out}\|_2$. Adding the two gives \eqref{eq:class_discrepancy}.
\end{proof}

\subsection{Proof of Proposition~\ref{prop:constant_integrability}}
\label{app:proof:constant_integrability}
\begin{proof}
For every $x\in\Xc$ and every $k,j$, we have
$|\phi_{k,j}(x)|\le |A_{k,j}^\top x+b_{k,j}|\le M\|A_{k,j}\|_1+|b_{k,j}|\le M\sqrt m\,\|A_{k,j}\|_2+|b_{k,j}|$.
Applying Minkowski's inequality to the averages gives
\begin{equation}
\label{eq:feature_envelope_bound}
\Bigl(\frac1N\sum_{j=1}^r\sum_{k=1}^N\sup_{x\in\Xc}|\phi_{k,j}(x)|^2\Bigr)^{1/2}
\le
M\sqrt m\,\Bigl(\frac1N\sum_{j=1}^r\sum_{k=1}^N\|A_{k,j}\|_2^2\Bigr)^{1/2}
+
\Bigl(\frac1N\sum_{j=1}^r\sum_{k=1}^N b_{k,j}^2\Bigr)^{1/2}.
\end{equation}
Let $c_{k,j}:=\|A_{k,j}\|_2+|b_{k,j}|$ and $S_N:=\sum_{j=1}^r\sum_{k=1}^N c_{k,j}$. Then
\begin{equation}
\Bigl(\frac1N\sum_{j=1}^r\sum_{k=1}^N\|A_{k,j}\|_2^2\Bigr)^{1/2}\le N^{-1/2}S_N,
\qquad
\Bigl(\frac1N\sum_{j=1}^r\sum_{k=1}^N b_{k,j}^2\Bigr)^{1/2}\le N^{-1/2}S_N.
\end{equation}
Together with \eqref{eq:class_discrepancy_constants}, these bounds give
\begin{equation}
\Lambda_W\le C_2S_N,
\qquad
C_2:=1+M\sqrt m,
\end{equation}
while
\begin{equation}
L_{\mathrm{class}}
\le
B_\star N^{-1/2}S_N.
\end{equation}
By Assumption~\ref{ass:feature_law}, every $c_{k,j}$ has the same finite expectation $\mu:=\Erf[c_{1,1}]$. Linearity gives $\Erf[S_N]=rN\mu$, irrespective of dependence between coordinate banks. Jensen's inequality gives
\begin{equation}
\label{eq:sum_log_bound}
\Erf\bigl[\log(1+S_N)\bigr]
\le
\log\bigl(1+rN\mu\bigr)
\le
\log N+\log(1+r\mu).
\end{equation}
For $N\ge1$, the preceding bounds and $N^{-1/2}\le1$ give
\begin{align}
\log(1+\Lambda_W)
&\le \log(1+C_2)+\log(1+S_N),\\
\log(1+L_{\mathrm{class}})
&\le \log(1+B_\star)+\log(1+S_N).
\end{align}
Taking expectations and applying \eqref{eq:sum_log_bound} proves
\eqref{eq:constant_log_growth} after enlarging $C_0$.
\end{proof}

\subsection{Proof of Theorem~\ref{thm:policy_error_propagation}}
\label{app:proof:policy_error_propagation}
\begin{proof}
Both systems use the same innovations and start from the same constant initial state. Their exogenous coordinates then coincide at every date, so the state difference is a difference of controlled coordinates. Both policies are evaluated on the same realization $\omega$, and the policies generate actions according to
\begin{equation}
    \alpha_t^{H_1}=H_1(X_t^{H_1},t),
    \qquad
    \alpha_t^{H_2}=H_2(X_t^{H_2},t).
\end{equation}
Define
\begin{equation}
    E_t:=\|X_t^{H_1}-X_t^{H_2}\|_1,
    \qquad
    A_t:=\|\alpha_t^{H_1}-\alpha_t^{H_2}\|_1.
\end{equation}
Adding and subtracting $H_1(X_t^{H_2},t)$ gives
\begin{align}
\label{eq:action_error_decomposition}
\|\alpha_t^{H_1}-\alpha_t^{H_2}\|_1
&=
\|H_1(X_t^{H_1},t)-H_2(X_t^{H_2},t)\|_1 \nonumber\\
&\leq
\|H_1(X_t^{H_1},t)-H_1(X_t^{H_2},t)\|_1 \nonumber\\
&\quad
+\|H_1(X_t^{H_2},t)-H_2(X_t^{H_2},t)\|_1 \nonumber\\
&\leq
L_1\|X_t^{H_1}-X_t^{H_2}\|_1+\varepsilon .
\end{align}
The Lipschitz property of the controlled-coordinate update, together with the equality of the exogenous coordinates, gives
\begin{align}
\label{eq:controlled_state_error_recursion}
\|X_{t+1}^{H_1}-X_{t+1}^{H_2}\|_1
&\leq
L_X\|X_t^{H_1}-X_t^{H_2}\|_1
+
L_A\|\alpha_t^{H_1}-\alpha_t^{H_2}\|_1 .
\end{align}
Equations~\eqref{eq:action_error_decomposition} and \eqref{eq:controlled_state_error_recursion} yield
\begin{equation}
\label{eq:action_error_before_state_recursion}
    A_t\leq L_1 E_t+\varepsilon
\end{equation}
and
\begin{equation}
\label{eq:state_error_recursion}
    E_{t+1}
    \leq
    \Gamma E_t+L_A\varepsilon .
\end{equation}
Since the initial states agree, $E_1=0$. Iterating the recursion gives
\begin{equation}
\label{eq:state_error_iterated}
    E_t
    \leq
    L_A\varepsilon\sum_{s=0}^{t-2}\Gamma^s,
    \qquad 2\le t\le T+1,
\end{equation}
With $E_1=0$, this proves the state bound through $T+1$. Substitution into $A_t\le L_1 E_t+\varepsilon$ proves the action bound through $T$.
\end{proof}

\subsection{Proof of Proposition~\ref{prop:final_loss_stability}}
\label{app:proof:final_loss_stability}
\begin{proof}
For a fixed $\beta\in\mathcal B$, define
\begin{equation}
    F_{H_1}(\beta)
    :=
    \Epath\left[\mathcal L(\beta,X^{H_1},\alpha^{H_1})\right].
\end{equation}
The trajectory-loss condition gives
\begin{align}
|F_{H_1}(\beta)-F_{H_2}(\beta)|
&\le \Epath\left[
|\mathcal L(\beta,X^{H_1},\alpha^{H_1})-\mathcal L(\beta,X^{H_2},\alpha^{H_2})|
\right]\nonumber\\
&\le K_\ell\left[
\sum_{t=1}^{T+1}\Epath\|X_t^{H_1}-X_t^{H_2}\|_1
+\sum_{t=1}^T\Epath\|\alpha_t^{H_1}-\alpha_t^{H_2}\|_1
\right].
\end{align}
The elementary bound
\begin{equation}
    \left|\inf_{\beta\in\mathcal B}F_{H_1}(\beta)
    -
    \inf_{\beta\in\mathcal B}F_{H_2}(\beta)\right|
    \leq
    \sup_{\beta\in\mathcal B}
    |F_{H_1}(\beta)-F_{H_2}(\beta)|
\end{equation}
then yields \eqref{eq:final_loss_stability}.
\end{proof}

\subsection{Proof of Corollary~\ref{cor:class_objective_stability}}
\label{app:proof:combined_objective_stability}
\begin{proof}
On a common realization of the innovations, write $d_t=\|X_t^{H_1}-X_t^{H_2}\|_1$ and $a_t=\|\alpha_t^{H_1}-\alpha_t^{H_2}\|_1$. Theorem~\ref{thm:policy_error_propagation} gives
\begin{equation}
d_t\le L_A\varepsilon\sum_{s=0}^{t-2}\Gamma^s,
\qquad
a_t\le\varepsilon+L_1L_A\varepsilon\sum_{s=0}^{t-2}\Gamma^s,
\end{equation}
with empty sums equal to zero. The trajectory-loss bound therefore gives, uniformly in $\beta\in\mathcal B$,
\begin{equation}
\begin{gathered}
|\mathcal L(\beta,X^{H_1},\alpha^{H_1})-\mathcal L(\beta,X^{H_2},\alpha^{H_2})|\\
\le K_\ell\left[\sum_{t=1}^{T+1}d_t+\sum_{t=1}^T a_t\right]\\
\le K_\ell\varepsilon\left[T+L_A(1+L_1)\sum_{s=0}^{T-1}(T-s)\Gamma^s\right]
=K_T(L_1)\varepsilon.
\end{gathered}
\end{equation}
Here we used $\sum_{t=1}^T\sum_{s=0}^{t-2}\Gamma^s\le\sum_{t=1}^{T+1}\sum_{s=0}^{t-2}\Gamma^s=\sum_{s=0}^{T-1}(T-s)\Gamma^s$. This proves \eqref{eq:combined_path_cost}. Taking expectations and using the infimum inequality proves $|J(H_1)-J(H_2)|\le K_T(L_1)\varepsilon$.
\end{proof}

\subsection{Proof of Theorem~\ref{thm:rfnn_objective_approx}}
\label{app:proof:rfnn_objective_approx}
\begin{proof}
First, the Fourier assumption in Theorem~\ref{thm:scalar_rfnn_approx} implies the benchmark-policy Lipschitz condition required by Theorem~\ref{thm:policy_error_propagation}. For each coordinate,
\begin{align}
    |H^\star_j(x)-H^\star_j(y)|
    &\leq
    \int_{\R^m}|e^{ix\cdot\xi}-e^{iy\cdot\xi}|\,|G_j(\xi)|d\xi \nonumber\\
    &\leq
    \|x-y\|_1
    \int_{\R^m}\|\xi\|_\infty |G_j(\xi)|d\xi .
\end{align}
The last integral is finite because Theorem~\ref{thm:scalar_rfnn_approx} assumes $\int_{\R^m}\max(1,\|\xi\|^2)|G_j(\xi)|d\xi<\infty$. Hence $H^\star$ is Lipschitz on $[-M,M]^m$ in the $\ell^1$ norm. One admissible choice is
\begin{equation}
\label{eq:reference_lipschitz_from_fourier}
    L_\star
    =
    \sum_{j=1}^r
    \int_{\R^m}\|\xi\|_\infty |G_j(\xi)|d\xi ,
\end{equation}
where $G_j$ is the Fourier density associated with $H^\star_j$. In particular, for any two policy inputs $x=(X,t)$ and $\widetilde x=(\widetilde X,t)$ in $\Xc$ with the same date, this gives $\|H^\star(X,t)-H^\star(\widetilde X,t)\|_1\le L_\star\|\widetilde X-X\|_1$, which is the Lipschitz hypothesis of Theorem~\ref{thm:policy_error_propagation} with $H_1=H^\star$ and $L_1=L_\star$.
Fix the RFNN feature draw and define the policy approximation error on the policy input space by
\begin{equation}
\label{eq:rfnn_policy_error_on_state_space}
    \Delta_N^\Theta
    :=
    \sup_{x\in\Xc}\|H_N^\Theta(x)-H^\star(x)\|_1.
\end{equation}
Apply Corollary~\ref{cor:class_objective_stability} with benchmark policy $H^\star$, comparison policy $H_N^\Theta$ and $\varepsilon=\Delta_N^\Theta$. Only the benchmark-policy Lipschitz constant $L_\star$ enters the resulting bound:
\begin{equation}
    |J(H_N^\Theta)-J(H^\star)|
    \leq
    K_T(L_\star)\Delta_N^\Theta .
\end{equation}
The expectation bound
\begin{equation}
    \Erf\left[
    \Delta_N^\Theta
    \right]
    \leq
    \frac{rC_{\max}}{\sqrt{N}}
\end{equation}
follows from Corollary~\ref{cor:vector_rfnn_approx}. The coordinatewise weight bound in that corollary also gives $H_N^\Theta\in\mathcal H_N$. Taking $\Erf$ expectations in the fixed-draw objective bound gives the claimed RFNN objective approximation rate.
\end{proof}

\subsection{Proof of Proposition~\ref{prop:param_loss_lipschitz}}
\label{app:proof:param_loss_lipschitz}
\begin{proof}
Fix a realization $\omega$ of the innovations and apply the pathwise estimate \eqref{eq:combined_path_cost} to $H_1=H_\theta$ and $H_2=H_{\theta'}$. The state-Lipschitz constant of $H_\theta$ is at most $L_{\rm class}$. With $\varepsilon:=\sup_{x\in\Xc}\|H_\theta(x)-H_{\theta'}(x)\|_1$, this gives, for every $\beta\in\mathcal B$,
\begin{equation}
|\ell(\omega;\theta,\beta)-\ell(\omega;\theta',\beta)|
\le K_T(L_{\rm class})\varepsilon.
\end{equation}
Substituting the two-term bound for $\varepsilon$ from \eqref{eq:class_discrepancy} proves the claim.
\end{proof}

\subsection{Proof of Theorem~\ref{thm:estimation_error}}
\label{app:proof:estimation_error}
\begin{proof}
For every $H$, the inequality
$|\inf_\beta u_\beta-\inf_\beta v_\beta|\le\sup_\beta|u_\beta-v_\beta|$
gives
\begin{equation}
\sup_{H\in\mathcal H_N}\bigl|\widehat J_n(H)-J(H)\bigr|
\le
\widehat Z_n,
\qquad
\widehat Z_n
:=
\sup_{\theta\in\mathcal W_N,\beta\in\mathcal B}
\Bigl|
\frac1n\sum_{i=1}^n\ell(\omega_i;\theta,\beta)
-
\E[\ell(\omega;\theta,\beta)]
\Bigr|.
\end{equation}
Let $\mathcal F:=\{\omega\mapsto\ell(\omega;\theta,\beta):
\theta\in\mathcal W_N,\beta\in\mathcal B\}$, fix $\ell_0\in\mathcal F$, and define
$\mathcal F_0:=\{\ell-\ell_0:\ell\in\mathcal F\}$. The compact parameter sets and parameter-Lipschitz bounds allow each supremum to be taken over a countable dense subset, so the suprema are measurable. The standard ghost-sample symmetrization argument used in the proof of Theorem~8 of \cite{BartlettMendelson2002}, together with $|\ell_0|\le U$, gives
\begin{equation}
\label{eq:Z_rademacher_reduction}
\E[\widehat Z_n]
\le
2\mathfrak R_n^{\mathrm{abs}}(\mathcal F_0)
+
\frac{2U}{\sqrt n},
\end{equation}
where
$\mathfrak R_n^{\mathrm{abs}}(\mathcal F_0)
:=
\E\sup_{f\in\mathcal F_0}
|n^{-1}\sum_{i=1}^n\varsigma_i f(\omega_i)|$.
The centered class satisfies $|f|\le2U$. It remains to bound its Rademacher complexity. With the empirical metric
$d_n(f,g):=(n^{-1}\sum_{i=1}^n|f(\omega_i)-g(\omega_i)|^2)^{1/2}$,
the process $X_f:=n^{-1}\sum_{i=1}^n\varsigma_i f(\omega_i)$ is sub-Gaussian for the metric $d_n/\sqrt n$ conditionally on the sample. Indeed, Hoeffding's lemma gives $\E_\varsigma\exp(\lambda(X_f-X_g))\le\exp(\lambda^2d_n(f,g)^2/(2n))$ for all $f,g$ and $\lambda$. Since $\mathcal F_0$ contains zero, $\sup_{f\in\mathcal F_0}|X_f|\le\sup_{f,g\in\mathcal F_0}(X_f-X_g)$, so Theorem~5.22 of \cite{Wainwright2019} with $\delta=0$ gives
\begin{equation}
\begin{aligned}
\mathfrak R_n^{\mathrm{abs}}(\mathcal F_0)
&\le
C\,\E\left[
\int_0^{2U/\sqrt n}
\sqrt{\log\mathcal N(\mathcal F_0,d_n/\sqrt n,\varepsilon)}\,d\varepsilon
\right]\\
&=
\frac{C}{\sqrt n}\,\E\left[
\int_0^{2U}
\sqrt{\log\mathcal N(\mathcal F_0,d_n,u)}\,du
\right].
\end{aligned}
\end{equation}
We now estimate the covering number. By Proposition~\ref{prop:param_loss_lipschitz} and the auxiliary-parameter Lipschitz bound in Assumption~\ref{ass:trajectory_loss},
\begin{equation}
|\ell(\omega;\theta,\beta)-\ell(\omega;\theta',\beta')|
\le
L_\ell^W\|W-W'\|_F
+
L_\ell^{\rm out}\|b^{\rm out}-b'^{\rm out}\|_2
+
K_\ell|\beta-\beta'|.
\end{equation}
Let $D_\beta:=\operatorname{diam}(\mathcal B)$. Rescale the parameter blocks according to
$\widetilde b^{\rm out}:=b^{\rm out}/B_{\rm out}$,
$\widetilde W:=\sqrt N W/B_\star$ when $B_\star>0$, and
$\widetilde\beta:=(\beta-\beta_0)/D_\beta$ when $D_\beta>0$, for some fixed $\beta_0\in\mathcal B$. If either radius is zero, the corresponding parameter block is fixed and is omitted from the following display and covering argument. Then
\begin{equation}
|\ell(\omega;\theta,\beta)-\ell(\omega;\theta',\beta')|
\le
\frac{B_\star L_\ell^W}{\sqrt N}
\|\widetilde W-\widetilde W'\|_F
+
B_{\rm out}L_\ell^{\rm out}
\|\widetilde b^{\rm out}-\widetilde b^{\rm out}{}'\|_2
+
K_\ell D_\beta
|\widetilde\beta-\widetilde\beta'|.
\end{equation}
Since $B_\star$, $B_{\rm out}$, $K_\ell$, and $D_\beta$ are fixed problem parameters, while
$\mathcal L_\ell=L_\ell^W+L_\ell^{\rm out}$, the loss is $C(1+\mathcal L_\ell)$-Lipschitz in the rescaled joint parameter. Since $\sum_{j=1}^r\|W_{j,\cdot}\|_2^2\le B_\star^2/N$, the rescaled parameter set is contained in a Euclidean ball of fixed radius in $\R^{r_N+1}$, and hence the standard volumetric bound gives
\begin{equation}
\mathcal N(\mathcal F,\|\cdot\|_\infty,u)
\le
\Bigl(
1+\frac{C(1+\mathcal L_\ell)}{u}
\Bigr)^{r_N+1}.
\end{equation}
Translation by $\ell_0$ does not change covering numbers, and
$d_n(f,g)\le\|f-g\|_\infty$, so
\begin{equation}
\mathcal N(\mathcal F_0,d_n,u)
\le
\Bigl(
1+\frac{C(1+\mathcal L_\ell)}{u}
\Bigr)^{r_N+1}.
\end{equation}
Substituting this estimate into the entropy integral yields
\begin{equation}
\mathfrak R_n^{\mathrm{abs}}(\mathcal F_0)
\le
\frac{C\sqrt{r_N+1}}{\sqrt n}
\int_0^{2U}
\sqrt{
\log\Bigl(
1+\frac{C(1+\mathcal L_\ell)}{u}
\Bigr)
}\,du.
\end{equation}
Using
$\int_0^{2U}\sqrt{\log(2U/u)}\,du=\sqrt{\pi}U$
and the inequality
$1+A/u\le(2U/u)(1+A/(2U))$ for $0<u\le2U$ gives
\begin{equation}
\mathfrak R_n^{\mathrm{abs}}(\mathcal F_0)
\le
\frac{CU\sqrt{r_N+1}}{\sqrt n}
\sqrt{\log\Bigl(2+\frac{\mathcal L_\ell}{U}\Bigr)}.
\end{equation}
Combining this with \eqref{eq:Z_rademacher_reduction} gives
\begin{equation}
\E[\widehat Z_n]
\le
\frac{CU\sqrt{r_N+1}}{\sqrt n}
\sqrt{\log\Bigl(2+\frac{\mathcal L_\ell}{U}\Bigr)}
+
\frac{2U}{\sqrt n}.
\end{equation}
The initial reduction to $\widehat Z_n$ proves \eqref{eq:estimation_error_bound}.
\end{proof}

\subsection{Proof of Proposition~\ref{prop:estimation_constant_growth}}
\label{app:proof:estimation_constant_growth}
\begin{proof}
Since $T$ is fixed, \eqref{eq:combined_cost_constant} gives
\begin{equation}
K_T(L_{\mathrm{class}})\le C(1+L_{\mathrm{class}})^T.
\end{equation}

By the definition of $\mathcal L_\ell$ in Theorem~\ref{thm:estimation_error}, with the fixed constant $\Lambda_{\rm out}$ absorbed into $C$,
\begin{equation}
1+\frac{\mathcal L_\ell}{U}
\le
C(1+\Lambda_W)(1+L_{\mathrm{class}})^T.
\end{equation}
Therefore
\begin{equation}
\log\Bigl(2+\frac{\mathcal L_\ell}{U}\Bigr)
\le
C+\log(1+\Lambda_W)+T\log(1+L_{\mathrm{class}}).
\end{equation}
Taking expectations and applying Proposition~\ref{prop:constant_integrability} gives \eqref{eq:loss_log_bound}. The definition of $\kappa_{\mathrm{est}}$ and Jensen's inequality then give \eqref{eq:kappa_growth}.
\end{proof}

\subsection[Proof of Corollary~\ref{cor:excess_risk}]{Proof of Corollary~\ref{cor:excess_risk}}
\label{app:proof:excess_risk}
\begin{proof}[Proof]
Fix the feature draw and write $Z_n:=\sup_{H\in\mathcal H_N}|\widehat J_n(H)-J(H)|$. Empirical optimality gives $\widehat J_n(\widehat H)\le\widehat J_n(H_N^\star)$, and $J(H)\le\widehat J_n(H)+|\widehat J_n(H)-J(H)|$ applied at $\widehat H$ and at $H_N^\star$ gives
\begin{equation}
\begin{aligned}
J(\widehat H)
&\le \widehat J_n(\widehat H)+\bigl|\widehat J_n(\widehat H)-J(\widehat H)\bigr|\\
&\le \widehat J_n(H_N^\star)+\bigl|\widehat J_n(\widehat H)-J(\widehat H)\bigr|\\
&\le J(H_N^\star)
+\bigl|\widehat J_n(H_N^\star)-J(H_N^\star)\bigr|
+\bigl|\widehat J_n(\widehat H)-J(\widehat H)\bigr|.
\end{aligned}
\end{equation}
This is \eqref{eq:excess_risk_decomposition}, and each deviation is at most $Z_n$. The feasible approximant $H_N^\Theta$ belongs to $\mathcal H_N$, so $J(H_N^\star)\le J(H_N^\Theta)$. Assumption~\ref{ass:trajectory_loss} gives $-U\le J(H)\le U$ for every policy under consideration, so $J(H_N^\star)$ and $J(H^\star)$ are finite. Theorem~\ref{thm:rfnn_objective_approx} therefore yields
\begin{equation}
\Erf[J(H_N^\star)-J(H^\star)]
\le \Erf\bigl[|J(H_N^\Theta)-J(H^\star)|\bigr]
\le \frac{K_T(L_\star) rC_{\max}}{\sqrt N}.
\end{equation}
Theorem~\ref{thm:estimation_error} and Proposition~\ref{prop:estimation_constant_growth} give
\begin{equation}
\Erf\,\E[Z_n]
\le \frac{\Erf[\kappa_{\mathrm{est}}]+2U}{\sqrt n}
\le C\sqrt{\frac{N(1+\log N)}{n}},
\end{equation}
where the last inequality uses $N(1+\log N)\ge1$ for $N\ge1$. Taking expectations in \eqref{eq:excess_risk_decomposition} and combining these estimates proves \eqref{eq:excess_risk_bound}.
\end{proof}

\section{Feature selection ablation and detailed PRFL setup}
\label{app:prfl_ablation}

PRFL trains only the output coefficients, and its optimizer settings follow Section~\ref{sec:numerical_training} except for the learning-rate schedule. In the Markovian examples the initial rate is held for the first $60\%$ of training, so the adaptive phase, which occupies the first half, runs entirely at that rate. We therefore run the first half as five PRFL cycles at the same rate, and the rest of the schedule is that of Section~\ref{sec:numerical_training}. All three methods use the weight decay $10^{-3}$ of that section. In Heston hedging no weight decay is used, and the adaptive phase extends beyond the first decay, so the schedule differs. The adaptive phase covers the first $800$ of the $1600$ epochs, in $16$ cycles of $50$ epochs. Within each cycle the learning rate is $0.005$ for $25$ epochs and $0.0005$ for the next $25$. After the adaptive phase the rate returns to $0.005$ for $200$ epochs and is then multiplied by $0.1$ every $200$ epochs, while the baselines multiply $0.005$ by $0.1$ every $400$ epochs.

Section~\ref{sec:prfl_benchmark_results} compares PRFL with RFNN and FTNN at matched retained width and training budget. PRFL is not trained in the same way as the baselines, since its dictionary grows and is pruned during the adaptive phase, so the comparison alone does not separate the effect of the selected features from the effect of the training procedure. To isolate the selection effect, we train an RFNN from scratch on the features PRFL retains and compare it with an RFNN of the same width trained on a uniformly random subset, with two independent draws per seed. Both policies use the same width, initialization, optimizer, learning-rate schedule, batch size and number of updates, and differ only in which features they use.

\begin{table}[!htbp]
\caption{Feature selection ablation. Mean objective over three seeds for PRFL, RFNN on a uniformly random subset of the retained width, and RFNN on the PRFL-selected features at the same width.}
\label{tab:prfl_ablation}
\centering
\small
\setlength{\tabcolsep}{6pt}
\begin{tabular}{lccc}
\hline
Setting & PRFL & Random & Selected \\
\hline
LQ $d=100$ & 5.4350 & 5.4493 & 5.4350 \\
Heston CVaR, $\lambda_{\mathrm{tc}}=0.005$ & 2.3345 & 2.4402 & 2.3341 \\
Heston square loss, $\lambda_{\mathrm{tc}}=0$ & 0.16691 & 0.17993 & 0.16537 \\
\hline
\end{tabular}
\end{table}

Table~\ref{tab:prfl_ablation} reports three representative settings. Both the selected features and PRFL outperform the random RFNN in every setting, and the selected features are lower than every random draw. The advantage over a random subset of the same width therefore comes from the selected features rather than the schedule.

\end{document}